\documentclass{article}

\usepackage{amssymb}
\usepackage{amsthm}
\usepackage[all]{xy}

\newtheorem{thm}{Theorem}[section]
\newtheorem{prop}[thm]{Proposition}
\newtheorem{lem}[thm]{Lemma}
\newtheorem{cor}[thm]{Corollary}
\newtheorem{ithm}{Theorem}
\newtheorem{iprop}{Proposition}

\theoremstyle{definition}
\newtheorem{dfn}[thm]{Definition}
\newtheorem{idfn}{Definition}

\theoremstyle{remark}
\newtheorem{rem}{Remark}
\newtheorem*{conventions}{Conventions}
\newtheorem*{acknowledgments}{Acknowledgments}

\newcommand{\C}{\mathbb{C}}
\newcommand{\R}{\mathbb{R}}
\newcommand{\T}{\mathbb{T}}
\newcommand{\Z}{\mathbb{Z}}

\newcommand{\g}{\mathfrak{g}}
\renewcommand{\AA}{\mathfrak{A}}

\newcommand{\F}{\mathcal{F}}
\renewcommand{\S}{\mathcal{S}}
\newcommand{\U}{\mathcal{U}}

\newcommand{\bF}{\bar{\F}}
\newcommand{\sC}{\mathcal{C}}

\newcommand{\Ad}{\mathrm{Ad}}
\newcommand{\Co}{\mathrm{Co}}
\newcommand{\id}{\mathrm{id}}

\newcommand{\Cech}{\v{C}ech {}}
\newcommand{\Poincare}{Poincar\'e {}}
\renewcommand{\i}{\sqrt{\! - \! 1}}
\renewcommand{\d}{\partial}
\renewcommand{\epsilon}{\varepsilon}
\def\u#1{ \underline{#1} }

\def\til#1{ \tilde{#1} }

\title{Equivariant smooth Deligne cohomology}

\author{Kiyonori Gomi
}

\date{}

\begin{document}

\maketitle

\begin{abstract}
On the basis of Brylinski's work, we introduce a notion of equivariant smooth Deligne cohomology group, which is a generalization of both the ordinary smooth Deligne cohomology and the ordinary equivariant cohomology. Using the cohomology group, we classify equivariant circle bundles with connection, and equivariant gerbes with connection.
\end{abstract}


\section{Introduction}

The \textit{smooth Deligne cohomology group} is a generalization of the cohomology group with coefficients $\Z$ in such a way that it incorporates the information of differential forms \cite{Bry1,De-F,E-V,Ga}. For a non-negative integer $N$ fixed, the smooth Deligne cohomology group $H^m(M, \F(N))$ of a smooth manifold $M$ is defined to be the hypercohomology of the complex of sheaves $\F(N)$ on $M$ given by
$$
\u{\T} \stackrel{\frac{1}{2\pi\i}d\log}{\longrightarrow}
\u{A}^1 \stackrel{d}{\longrightarrow} 
\u{A}^2 \stackrel{d}{\longrightarrow} 
\cdots \stackrel{d}{\longrightarrow}
\u{A}^N \longrightarrow
0 \longrightarrow \cdots,
$$
where $\u{\T}$ is the sheaf of germs of smooth functions with values in the unit circle $\T = \{ z \in \C |\ |z| = 1 \}$, and $\u{A}^q$ is the sheaf of germs of smooth differential $q$-forms with values in $\R$. 

The usefulness of the smooth Deligne cohomology groups would be best understood by their geometric interpretations. Recall the following geometric interpretations of ordinary cohomology groups of low degree.

\begin{iprop} \label{iprop:oridinary_coh}
Let $M$ be a smooth manifold.

(a)(Kostant \cite{Ko}, Weil \cite{We}) The isomorphism classes of principal $\T$-bundles (Hermitian line bundles) over $M$ are classified by $H^1(M, \u{\T}) \cong H^2(M, \Z)$.

(b)(Giraud \cite{Gi}) The isomorphism classes of gerbes over $M$ are classified by $H^2(M, \u{\T}) \cong H^3(M, \Z)$.
\end{iprop}

In the above, a ``gerbe'' means a \textit{gerbe with band $\u{\T}$} \cite{Bry1,Gi}. We remark that a class in $H^2(M, \u{\T}) \cong H^3(M, \Z)$ admits the other geometric interpretations. (See \cite{Ga,H,Mu} for example.) We also remark that an interpretation of cohomology groups of the degree equal or greater than four requires notions of higher gerbes \cite{Bry-M,C-M-W,F2,St}.

The geometric interpretation of smooth Deligne cohomology groups is obtained as a generalization of Proposition \ref{iprop:oridinary_coh}.

\begin{iprop}[Brylinski \cite{Bry1}] \label{iprop:Deligne_coh}
Let $M$ be a smooth manifold.

(a) The isomorphism classes of principal $\T$-bundles with connection over $M$ are classified by $H^1(M, \F(1))$.

(b) The isomorphism classes of gerbes with connective structure and curving over $M$ are classified by $H^2(M, \F(2))$.
\end{iprop}

When a Lie group $G$ acts on a smooth manifold $M$, one often consider equivariant cohomology groups to include the information of the group action. The standard definition of equivariant cohomology groups is the Borel construction \cite{A-B}, that is, $H^m_G(M, \Z) = H^m(EG \times_G M, \Z)$. However, to obtain an equivariant generalization of Proposition \ref{iprop:oridinary_coh}, we consider the other cohomology.

In the case that a Lie group $G$ acts on a smooth manifold $M$, we have a simplicial manifold (\cite{Du,Se}) $G^{\bullet} \times M = \{ G^p \times M \}_{p \ge 0}$. The family of sheaves $\{ \u{\T}_{G^p \times M} \}_{p \ge 0}$ gives rise to a \textit{simplicial sheaf} \cite{De} over $G^{\bullet} \times M$, where $\u{\T}_{G^p \times M}$ is the sheaf of germs of $\T$-valued smooth functions on $G^p \times M$. We denote the hypercohomology of this simplicial sheaf by $H^m(G^{\bullet} \times M, \u{\T})$. Note that, if $G$ is compact, then there is an isomorphism $H^m(G^{\bullet} \times M, \u{\T}) \cong H^{m+1}_G(M, \Z)$ for $m \ge 1$. 

\begin{iprop}[Brylinski \cite{Bry2}] \label{iprop:equiv_coh}
Let $G$ be a Lie group acting on a smooth manifold $M$.

(a) The isomorphism classes of $G$-equivariant principal $\T$-bundles over $M$ are classified by $H^1(G^{\bullet} \times M, \u{\T})$.

(b) The isomorphism classes of $G$-equivariant gerbes over $M$ are classified by $H^2(G^{\bullet} \times M, \u{\T})$.
\end{iprop}

The purpose of the present paper is to formulate ``equivariant smooth Deligne cohomology groups'' which allow one to have an equivariant generalization of Proposition \ref{iprop:Deligne_coh}. The formulation owes the fundamental ideas to Brylinski's paper \cite{Bry2}. Let $G^{\bullet} \times M$ be the simplicial manifold associated to the group action of $G$ on $M$. On each piece $G^p \times M$ of this simplicial manifold, we introduce a complex of sheaves $\bF(N)_{G^p \times M}$ by
$$
\u{\T} \stackrel{\frac{1}{2\pi\i}d\log}{\longrightarrow}
\u{A}^1_{rel} \stackrel{d}{\longrightarrow} 
\u{A}^2_{rel} \stackrel{d}{\longrightarrow} 
\cdots \stackrel{d}{\longrightarrow}
\u{A}^N_{rel} \longrightarrow
0 \longrightarrow \cdots.
$$
In the above, we denote by $\u{A}^q_{rel}$ the sheaf of germs of relative differential $q$-forms on $G^p \times M$ with respect to the fibration $\pi: G^p \times M \to G^p \times pt$, where $pt$ is the manifold consisting of a single point. The family $\{ \bF(N)_{G^p \times M} \}_{p \ge 0}$ gives rise to a complex of simplicial sheaves $\bF(N)$ on $G^{\bullet} \times M$. 

\begin{idfn}
Let $G$ be a Lie group acting on a smooth manifold $M$. We define the \textit{equivariant smooth Deligne cohomology group} $H^m(G^{\bullet} \times M, \bF(N))$ to be the hypercohomology of the complex of simplicial sheaves $\bF(N)$ on $G^{\bullet} \times M$.
\end{idfn}

When $G$ is a finite group, the equivariant smooth Deligne cohomology group above coincides with the \textit{Deligne cohomology group for the orbifold $M/G$} introduced by Lupercio and Uribe \cite{L-U}.

Now we state the generalization of Proposition \ref{iprop:Deligne_coh} and Proposition \ref{iprop:equiv_coh}.

\begin{ithm}
Let $G$ be a Lie group acting on a smooth manifold $M$.

(a) The isomorphism classes of $G$-equivariant principal $\T$-bundles with $G$-invariant connection over $M$ are classified by $H^1(G^{\bullet} \times M, \bF(1))$.

(b) The isomorphism classes of $G$-equivariant gerbes with $G$-invariant connective structure and $G$-invariant curving over $M$ are classified by $H^2(G^{\bullet} \times M, \bF(2))$.
\end{ithm}

By virtue of this classification, we can extract various informations of equivariant gerbes with connection (i.e.\ $G$-invariant connective structure and $G$-invariant curving) from the knowledge of cohomology groups. For example, we obtain obstruction classes for an ordinary gerbe with connection to being equivariant (Corollary \ref{cor:obstruction_moduli_gerbe_connection}). In \cite{Go}, a similar cohomological method is applied to the study of relationships between $G$-equivariant gerbes with connection over $M$ and gerbes with connection over the quotient space $M/G$.

\bigskip

The organization of this paper is as follows. 

In Section \ref{sec:SDC}, we define the smooth Deligne cohomology group, and review its basic properties. 

In Section \ref{sec:EDC}, we introduce the equivariant smooth Deligne cohomology group $H^m(G^{\bullet} \times M, \bF(N))$. For this aim, we define the simplicial manifold $G^{\bullet} \times M$ and explain the hypercohomology of a simplicial sheaf. The \Cech cohomology description is also explained. 

In Section \ref{sec:properties_EDC}, we study the equivariant smooth Deligne cohomology group in the relations with the ordinary smooth Deligne cohomology, the ordinary equivariant cohomology, and the group of invariant differential forms. In addition, we construct a homomorphism to the equivariant de Rham cohomology. 

In Section \ref{sec:equiv_T_bundle_and_gerbe}, we first give the classification of equivariant principal $\T$-bundles with connection. For this aim, we give the ``simplicial'' formulation of equivariant principal $\T$-bundles, following the paper \cite{Bry2}. We next give the notion of equivariant gerbes with connection. (To save the pages, we drop the definition of gerbe itself. We refer the reader to \cite{Bry1,Bry-M}.) Then we state the classification of equivariant gerbes with connection by using the equivariant smooth Deligne cohomology. Some results on equivariant gerbes with connection are derived as a simple application of results in Section \ref{sec:properties_EDC}.

\begin{conventions}
Throughout this paper, we make a convention that a ``smooth manifold'' means a paracompact smooth manifold modeled on a topological vector space which is Hausdorff and locally convex. We also assume the existence of a partition of unity. Examples of such a manifold cover not only all the finite dimensional smooth manifolds, but also a sort of infinite dimensional manifolds. (The most important example of the infinite dimensional case would be the loop space of a finite dimensional smooth manifold. See \cite{Bry1} for detail.) 

We also make a convention that a ``Lie group'' means a Lie group whose underlying smooth manifold is of the type above. When a Lie group $G$ acts on a smooth manifold $M$, we assume that the action is smooth. We denote the action by juxtaposition: we write $gx \in M$ for $g \in G$ and $x \in M$. We always denote by $e \in G$ the unit element of the Lie group.
\end{conventions}


\section{Review of Smooth Deligne cohomology}
\label{sec:SDC}

This section is devoted to recalling ordinary smooth Deligne cohomology groups \cite{Bry1,De-F,E-V,Ga}. 

\subsection{Smooth Deligne cohomology groups}

Let $M$ be a smooth manifold. We denote by $\u{\T}_M$ the sheaf of germs of smooth functions with values in $\T = \{ u \in \C |\ |u| = 1 \}$. For a non-negative integer $q$, we denote by $\u{A}^q_M$ the sheaf of germs of $\R$-valued smooth $q$-forms on $M$. 

\begin{dfn}[\cite{Bry1}]
Let $N$ be a non-negative integer. 

(a) We define the \textit{smooth Deligne complex} $\F(N)_M$ to be the following complex of sheaves on $M$:
$$
\F(N)_M : \
\u{\T}_M \stackrel{\frac{1}{2\pi\i}d\log}{\longrightarrow}
\u{A}^1_M \stackrel{d}{\longrightarrow} 
\u{A}^2_M \stackrel{d}{\longrightarrow} 
\cdots \stackrel{d}{\longrightarrow}
\u{A}^N_M \longrightarrow
0 \longrightarrow \cdots,
$$
where the sheaf $\u{\T}_M$ is located at degree 0 in the complex. 

(b) The \textit{smooth Deligne cohomology group} $H^p(M, \F(N)_M)$ is defined to be the hypercohomology group of the smooth Deligne complex.
\end{dfn}

We often omit the subscripts of $\u{\T}_M, \u{A}^q_M$ and $\F(N)_M$.

\begin{rem}
Let $\Z(N)_D^\infty$ be a complex of sheaves given by
$$
\Z(N)_D^\infty : \
\Z \stackrel{i}{\longrightarrow}
\u{A}^0 \stackrel{d}{\longrightarrow} 
\u{A}^1 \stackrel{d}{\longrightarrow} 
\u{A}^2 \stackrel{d}{\longrightarrow}
\cdots \stackrel{d}{\longrightarrow}
\u{A}^{N-1} \longrightarrow
0 \longrightarrow \cdots,
$$
where we regard $\Z$ as a constant sheaf over $M$. The smooth Deligne cohomology often refers to the hypercohomology $H^p(M, \Z(N)_D^\infty)$. Since $\Z(N)_D^\infty$ is quasi-isomorphic to $\F(N-1)$ under a shift of degree, we have $H^p(M, \Z(N)_D^\infty) \cong H^{p-1}(M, \F(N-1))$.
\end{rem}

The smooth Deligne complex fits into the following short exact sequences of complexes of sheaves on $M$:
\begin{eqnarray}
& 0 \to
\{ \u{\T} \to \u{A}^1 \to \cdots \to \u{A}^N_{cl} \} \to
\F(N) \stackrel{d}{\to}
\{ 0 \to \cdots \to 0 \to \u{A}^{N+1}_{cl} \} \to 0, &
\label{exact_seq:Deligne_to_forms:manifold} \\
& 0 \to
\{ 0 \to \u{A}^1 \to \cdots \to \u{A}^N \} \to
\F(N) \to
\{ \u{\T} \to 0 \to \cdots \to 0 \} \to 0, &
\label{exact_seq:Deligne_to_integral:manifold}
\end{eqnarray}
where $\u{A}^q_{cl}$ is the sheaf of germs of closed $q$-forms on $M$. By the \Poincare lemma \cite{B-T,Bry1}, there exists a quasi-isomorphism
$$
\{ \T \to 0 \to \cdots \to 0 \} \to
\{ \u{\T} \to \u{A}^1 \to \cdots \to \u{A}^N_{cl} \},
$$
where $\T$ means the constant sheaf on $M$.

\begin{prop}[\cite{Bry1}] \label{prop:Deligne_coh:manifold}
Let $N$ be a positive integer.

(a) If $0 \le p < N$, then $H^p(M, \F(N))$ is isomorphic to $H^p(M, \T)$.

(b) If $p = N$, then $H^N(M, \F(N))$ fits into the following exact sequences:
\begin{eqnarray*}
& 0 \longrightarrow
H^N(M, \T) \longrightarrow
H^N(M, \F(N)) \stackrel{d}{\longrightarrow}
A^{N+1}(M)_0 \longrightarrow 0, & \\
& 0 \longrightarrow
A^N(M) / A^N(M)_0 \longrightarrow
H^N(M, \F(N)) \longrightarrow
H^{N+1}(M, \Z) \longrightarrow 0, &
\end{eqnarray*}
where $A^q(M)_0$ is the group of closed integral $q$-forms on $M$. 

(c) If $N < p$, then $H^p(M, \F(N))$ is isomorphic to $H^p(M, \u{\T}) \cong H^{p+1}(M, \Z)$.
\end{prop}

\begin{proof}
By the \Poincare lemma, we can take $(\u{A}^{N+1+*}, d)$ as a resolution of $\u{A}^{N+1}_{cl}$. Since $M$ is assumed to admit a partition unity, $\u{A}^q$ is \textit{soft} \cite{Bry1} and we have
$$
H^p(M, 0 \to \cdots \to 0 \to \u{A}^{N+1}_{cl}) =
\left\{
\begin{array}{cl}
0, & (0 \le p < N), \\
A^{N+1}(M)_{cl}, & (p = N), \\
H^{p+1}(M, \R), & (N < p). \\
\end{array}
\right.
$$
By the long exact sequence associated with (\ref{exact_seq:Deligne_to_forms:manifold}), we have the isomorphism in (a) and the first exact sequence in (b). Similarly, by a computation, we obtain
$$
H^p(M, 0 \to \u{A}^1 \to \cdots \to \u{A}^N) =
\left\{
\begin{array}{cl}
A^N(M)/ dA^{N-1}(M), & ( p = N), \\
0, & (N < p), \\
\end{array}
\right.
$$
where $d A^{N-1}(M)$ is the image of $d : A^{N-1}(M) \to A^{N}(M)$. Thus, the long exact sequence associated with (\ref{exact_seq:Deligne_to_integral:manifold}) gives the second exact sequence in (b) and the isomorphism in (c).
\end{proof}

It is easy to see that $H^0(M, \F(0)) = H^0(M, \u{\T})$ is identified with $C^{\infty}(M, \T)$, the group of $\T$-valued smooth functions on $M$.

\subsection{Relative Deligne cohomology groups}

Let $\pi : E \to B$ be a smooth fibration with fiber $F$. The sheaf $\u{A}^p_B$ of germs of $p$-forms on $B$ is naturally a module over $\u{\R}_B = \u{A}^0_B$, so that the inverse image sheaf $\pi^{-1}\u{A}^p_B$ is a $\pi^{-1}\u{\R}_B$-module. For a positive integer $p$, we define a subsheaf $F^p\!\u{A}^q_E$ of $\u{A}^q_E$ by setting $F^p\!\u{A}^q_E = \pi^{-1}\u{A}^p_B \otimes_{\pi^{-1}\u{\R}_B} \u{A}^{q-p}_E$, where $\u{A}^{q-p}_E$ is regarded as a $\pi^{-1}\u{\R}_B$-module through the natural homomorphism $\pi^{-1}\u{\R}_B \to \u{\R}_E$.

For an open subset $U \subset E$, the group $F^p\!\u{A}^q_E(U)$ consists of those $q$-forms $\omega$ on $U$ satisfying $\iota_{V_1} \cdots \iota_{V_{q-p+1}} \omega = 0$ for tangent vectors $V_1, \ldots, V_{q-p+1}$ at $x \in U$ such that $\pi_*V_j = 0$. If $\{ x_i \}$ and $\{ y_j \}$ are systems of local coordinates of $B$ and $F$ respectively, then the $q$-form $\omega$ has a local expression
$$
\omega = 
\sum_{r \ge p} 
\sum_{{i_1, \ldots, i_r}, \atop {j_1, \ldots, j_{q-r}}}
f_{I, J}(x, y) 
dx_{i_1} \wedge \cdots \wedge dx_{i_r} \wedge
dy_{j_1} \wedge \cdots dy_{j_{q-r}}.
$$

As is clear, we have a filtration $\u{A}^q_E \supset F^1\!\u{A}^q_E \supset F^2\!\u{A}^q_E \supset \cdots \supset F^q\!\u{A}^q_E \supset 0$. Thus, the smooth Deligne complex $\F(N) = \F(N)_E$ on $E$ admits a filtration 
$$
\F(N) \supset
F^1\!\F(N) \supset
F^2\!\F(N) \supset 
\cdots \supset
F^N\!\F(N) \supset 0
$$
associated to the fibration $\pi : E \to B$. 
(See Figure \ref{fig:subcomplex_Deligne}.)

\begin{figure}[htbp]
$$
\xymatrix@C=15pt@R=2pt{
\u{\T} \ar[r] &
\u{A}^1 \ar[r] &
\u{A}^2 \ar[r] &
\cdots \ar[r] &
\u{A}^{N-1} \ar[r] &
\u{A}^N \ar[r] &
0 \ar[r] &
\cdots, \\
\cup &
\cup &
\cup &
\cdots &
\cup &
\cup &
\cup & \\
0 \ar[r] &
F^1\!\u{A}^1 \ar[r] &
F^1\!\u{A}^2 \ar[r] &
\cdots \ar[r] &
F^1\!\u{A}^{N-1} \ar[r] &
F^1\!\u{A}^N \ar[r] &
0 \ar[r] &
\cdots, \\
\cup &
\cup &
\cup &
\cdots &
\cup &
\cup &
\cup & \\
0 \ar[r] &
0 \ar[r] &
F^2\!\u{A}^2 \ar[r] &
\cdots \ar[r] &
F^2\!\u{A}^{N-1} \ar[r] &
F^2\!\u{A}^N \ar[r] &
0 \ar[r] &
\cdots, \\
\cup &
\cup &
\cup &
\cdots &
\cup &
\cup &
\cup & \\
\vdots &
\vdots &
\vdots &
\cdots &
\vdots &
\vdots &
\vdots &
\cdots, \\
\cup &
\cup &
\cup &
\cdots &
\cup &
\cup &
\cup & \\
0 \ar[r] &
0 \ar[r] &
0 \ar[r] &
\cdots \ar[r] &
0 \ar[r] &
F^N\!\u{A}^N \ar[r] &
0 \ar[r] &
\cdots, \\
\cup &
\cup &
\cup &
\cdots &
\cup &
\cup &
\cup & \\
0 \ar[r] &
0 \ar[r] &
0 \ar[r] &
\cdots \ar[r] &
0 \ar[r] &
0 \ar[r] &
0 \ar[r] &
\cdots.
}
$$
\caption{The filtration of $\F(N)_E$}
\label{fig:subcomplex_Deligne}
\end{figure}

\smallskip

The \textit{relative Deligne complex} $\bF(N)$ is defined by $\bF(N) = \F(N)/ F^1\!\F(N)$. When we denote by $\u{A}^q_{rel} = \u{A}^q/F^1\!\u{A}^q$ the sheaf of germs of \textit{relative differential $q$-forms} with respect to the fibration $\pi : E \to B$, the relative Deligne complex $\bF(N)$ is expressed as
$$
\bF(N) : \
\u{\T} \stackrel{\frac{1}{2\pi\i}d\log}{\longrightarrow}
\u{A}^1_{rel} \stackrel{d}{\longrightarrow} 
\u{A}^2_{rel} \stackrel{d}{\longrightarrow} 
\cdots \stackrel{d}{\longrightarrow}
\u{A}^N_{rel} \longrightarrow
0 \longrightarrow \cdots.
$$

We call the hypercohomology group $H^m(E, \bF(N))$ the \textit{relative Deligne cohomology} \cite{Bry1,Bry-M}. Clearly, if $B$ consists of a single point, then $\bF(N) = \F(N)$, so that $H^m(E, \bF(N)) = H^m(E, \F(N))$.

The relative Deligne complex fits into the following short exact sequences of complexes of sheaves on $E$ similar to (\ref{exact_seq:Deligne_to_forms:manifold}) and (\ref{exact_seq:Deligne_to_integral:manifold}):
\begin{eqnarray}
&
\xymatrix@C=7pt{
0 \ar[r] &
\{ \u{\T} \to \u{A}^1_{rel} \to \cdots \to \u{A}^N_{rel, cl} \} \ar[r] &
\bF(N) \ar[r]^-{d} &
\{ 0 \to \cdots \to 0 \to \u{A}^{N+1}_{rel, cl} \} \ar[r] & 0,
} 
& \label{exact_seq:Deligne_to_forms:relative} \\
&
0 \to
\{ 0 \to \u{A}^1_{rel} \to \cdots \to \u{A}^N_{rel} \} \to
\bF(N) \to
\{ \u{\T} \to 0 \to \cdots \to 0 \} \to 0,
& \label{exact_seq:Deligne_to_integral:relative}
\end{eqnarray}
where $\u{A}^q_{rel,cl}$ is the sheaf of germs of closed relative $q$-forms on $M$. By the relative version of the \Poincare lemma \cite{Bry1}, there exists a quasi-isomorphism
$$
\{ \pi^{-1}\u{\T} \to 0 \to \cdots \to 0 \} \to
\{ \u{\T} \to \u{A}^1_{rel} \to \cdots \to \u{A}^N_{rel,cl} \},
$$
where $\pi^{-1}\u{\T} = \pi^{-1}\u{\T}_B$ is the inverse image of the sheaf $\u{\T}_B$ under $\pi : E \to B$.

We can obtain the following properties of the relative Deligne cohomology by the same method as that used in the proof of Proposition \ref{prop:Deligne_coh:manifold}.

\begin{prop}
Let $N$ be a positive integer.

(a) If $0 \le p < N$, then $H^p(E, \bF(N))$ is isomorphic to $H^p(E, \pi^{-1}\u{\T})$.

(b) If $p = N$, then $H^N(E, \bF(N))$ fits into the exact sequences:
$$
0 \to H^N(E, \pi^{-1}\u{\T}) \to H^N(E, \bF(N)) \to A^{N+1}(E)_{rel, cl} \to 
H^{N+1}(E, \pi^{-1}\u{\T}),
$$
$$
H^N(E, \Z) \to A^N(E)_{rel} / d A^{N-1}(E)_{rel} \to
H^N(E, \bF(N)) \to 
H^{N+1}(E, \Z) \to 0.
$$

(c) If $N < p$, then $H^p(E, \bF(N))$ is isomorphic to  $H^p(E, \u{\T}) \cong H^{p+1}(E, \Z)$.
\end{prop}

For computations of $H^m(E, \pi^{-1}\u{\T}_B)$, the next lemma is useful.

\begin{lem}
Let $\pi : E \to B$ be a smooth fibration with fiber $F$. There exists a spectral sequence converging to a graded quotient of $H^m(E, \pi^{-1}\u{\T}_B)$ with its $E_2$-term given by
$$
E_2^{p, q} = H^p(B, \u{\T}_B \times_{\T} \u{H}^q(F, \T)),
$$
where $\u{H}^q(F, \T)$ is the sheaf on $B$ associated with the presheaf given by the assignment of $H^q(\pi^{-1}(V), \T)$ to an open set $V \subset B$.
\end{lem}

\begin{proof}
The proof is essentially the same as that of a lemma in \cite{Bry1} (1.6.9 Lemma, p.\ 59). By virtue of the fibration $\pi : E \to B$, we have the Leray spectral sequence \cite{Bry1}, that is, a spectral sequence converging to a graded quotient of $H^m(E, \pi^{-1}\u{\T})$ with its $E_2$-term given by $E_2^{p, q} = H^p(B, \u{\mathcal{H}}^q)$. Here $\u{\mathcal{H}}^q$ is the sheaf associated with the presheaf given by the assignment to an open set $V \subset B$ of the group $H^q(\pi^{-1}(V), \pi^{-1}\u{\T})$. Now we prove that the sheaf $\u{\T} \times_{\T} \u{H}^q(F, \T)$ is naturally isomorphic to $\u{\mathcal{H}}^q$. For this aim, it suffices to show that $\u{\T}(V) \times_{\T} H^q(\pi^{-1}(V), \T)$ is naturally isomorphic to $H^q(\pi^{-1}(V), \pi^{-1}\u{\T})$ for a sufficiently small contractible open set $V \subset B$. Let $\{ U_\alpha \}_{\alpha \in \AA}$ be a good cover of $F$. If we take and fix a local trivialization $\pi^{-1}(V) \simeq V \times F$, then $\{ V \times U_\alpha\}_{\alpha \in \AA}$ is a good cover of $\pi^{-1}(V)$. We can see that $\u{\T}(V) \times_{\T} H^q(\pi^{-1}(V), \T)$ is the $q$th cohomology of the \Cech complex $\prod_{\alpha_0, \ldots, \alpha_j} \Gamma(V, \u{\T}) \times_\T \Gamma(U_{\alpha_0} \cap \cdots \cap U_{\alpha_j}, \T)$. Because $H^q(\pi^{-1}(V), \pi^{-1}\u{\T})$ is computed from the same complex, we obtain the natural isomorphism.
\end{proof}


\section{Equivariant smooth Deligne cohomology}
\label{sec:EDC}

We formulate equivariant smooth Deligne cohomology groups here. As is mentioned, the formulation owes the basic ideas to Brylinski's paper \cite{Bry2}.

\subsection{Simplicial manifolds associated to group actions}

First of all, we introduce a certain \textit{simplicial manifold} \cite{Se} associated to a manifold with a group action. 

Let $G$ be a Lie group acting on a smooth manifold $M$ by left. Then we have a simplicial manifold $G^{\bullet} \times M = \{ G^p \times M \}_{p \ge 0}$, where the face maps $\d_i : G^{p+1} \times M \to G^{p} \times M$, $(i = 0, \ldots p+1)$ are given by
$$
\d_i(g_1, \ldots, g_{p+1}, x) =
\left\{
\begin{array}{ll}
(g_2, \ldots, g_{p+1}, x), & i = 0 \\
(g_1, \ldots, g_{i-1}, g_i g_{i+1}, g_{i+2}, \ldots, g_{p+1}, x), & i = 1, \ldots, p \\
(g_1, \ldots, g_p, g_{p+1} x), & i = p + 1,
\end{array}
\right.
$$
and the degeneracy maps $s_i : G^p \times M \to G^{p+1} \times M$, $(i = 0, \ldots p)$ by
$$
s_i(g_1, \ldots, g_p, x) = (g_1, \ldots, g_i, e, g_{i+1}, \ldots, g_p, x).
$$
These maps obey the following relations:
\begin{eqnarray}
\d_i \circ \d_j & = & \d_{j-1} \circ \d_i, \quad (i < j), \label{rel1} \\
s_i \circ s_j & = & s_{j+1} \circ s_i, \quad (i \le j), \label{rel2} \\
\d_i \circ s_j & = &
\left\{
\begin{array}{ll}
s_{j-1} \circ \d_i, & \quad (i < j), \\
\id, & \quad (i = j, j+1), \\
s_j \circ \d_{i-1}, &\quad (i > j+1).
\end{array}
\right. \label{rel3}
\end{eqnarray}

To a simplicial manifold, we can associate a topological space called the \textit{realization} \cite{Du,Ga,Se}. The realization of $G^{\bullet} \times M$ is identified with the homotopy quotient (\cite{A-B}): $|G^{\bullet} \times M| \cong EG \times_G M$, where $EG$ is the total space of the universal bundle for $G$. This can be seen by the fact that $EG$ is obtained as the realization of $G^{\bullet} \times G$, where $G$ acts on itself by the left translation. 

Note that the classifying space $BG$ is also obtained as the realization of $G^{\bullet} \times pt$, where $pt$ is the space consisting of a single point on which $G$ acts trivially. We denote by $\pi : G^{\bullet} \times M \to G^{\bullet} \times pt$ the map of simplicial manifolds given by the projection $\pi : G^p \times M \to G^p \times pt$.

\subsection{Definition of equivariant smooth Deligne cohomology}

In order to define equivariant smooth Deligne cohomology groups, we briefly explain the notion of a \textit{sheaf on a simplicial manifold} (a \textit{simplicial sheaf}, for short) and its cohomology group \cite{De}. 

Let $G^{\bullet} \times M$ be the simplicial manifold associated to an action of a Lie group $G$ on a smooth manifold $M$. We define a \textit{simplicial sheaf} on $G^{\bullet} \times M$ to be a family of sheaves $\S^{\bullet} = \{ \S^p \}_{p \ge 0}$, where $\S^p$ is a sheaf on $G^p \times M$ such that there are homomorphisms $\til{\d}_i : \d_i^{-1} \S^p \to \S^{p+1}$ and $\til{s}_i : s_i^{-1} \S^{p+1} \to \S^p$ obeying the same relations as (\ref{rel1}), (\ref{rel2}) and (\ref{rel3}).

\smallskip

For each $p$, let $I^{p, *}$ be an injective resolution of the sheaf $\S^{p}$ on $G^p \times M$. 
$$
\S^{p} 
\hookrightarrow
I^{p, 0} \stackrel{\delta}{\longrightarrow}
I^{p, 1} \stackrel{\delta}{\longrightarrow}
I^{p, 2} \stackrel{\delta}{\longrightarrow}
\cdots. 
$$
We denote this by $I^{*, *}$, and call an injective resolution of the simplicial sheaf $\S^{\bullet}$. The homomorphism $\til{\d}_i : \d_i^{-1} \S^p \to \S^{p+1}$ induces a homomorphism $\d_i^* : \Gamma(G^p \times M, I^{p,q}) \to \Gamma(G^{p+1} \times M, I^{p+1,q})$. Combining these homomorphisms, we define a homomorphism $\d : \Gamma(G^p \times M, I^{p,q}) \to \Gamma(G^{p+1} \times M, I^{p+1,q})$ by $\d = \sum_{j=0}^{p+1} (-1)^j\d_j^*$. This homomorphism satisfies $\d \circ \d = 0$, because of (\ref{rel1}).

Now the hypercohomology of the simplicial sheaf $\S^{\bullet}$ is defined to be the cohomology of the double complex $(\Gamma(G^i \times M, I^{i, j}), \d, \delta)$. This cohomology group is independent of the choice of an injective resolution $I^{*, *}$. The independence is shown by the same method as in the case of the ordinary hypercohomology. We denote the cohomology by $H^*(G^{\bullet} \times M, \S^{\bullet})$. 

The notion of a complex of simplicial sheaves and of its hypercohomology group are defined in a similar fashion.

\begin{dfn}
Let $G$ be a Lie group acting on a smooth manifold $M$. 

(a) We define a complex of simplicial sheaves $\F(N)_{G^{\bullet} \times M}$ on $G^{\bullet} \times M$ by the family of the smooth Deligne complex $\{ \F(N)_{G^p \times M} \}_{p \ge 0}$. The homomorphisms $\til{\d}_i : \d_i^{-1} \F(N)_{G^p \times M} \to \F(N)_{G^{p+1} \times M}$ and $\til{s}_i : s_i^{-1} \F(N)_{G^{p+1} \times M} \to \F(N)_{G^p \times M}$ are the natural ones.

(b) We define a subcomplex $F^1\!\F(N)_{G^{\bullet} \times M}$ of $\F(N)_{G^{\bullet} \times M}$ by the family $\{ F^1\!\F(N)_{G^p \times M} \}_{p \ge 0}$, where $F^1\!\F(N)_{G^p \times M}$ is the subcomplex of $\F(N)_{G^p \times M}$ associated to the fibration $\pi : G^p \times M \to G^p \times pt$.

(c) We define a complex of simplicial sheaves $\bF(N)_{G^{\bullet} \times M}$ on $G^{\bullet} \times M$ by the family $\{ \bF(N)_{G^p \times M} \}_{p \ge 0}$, where $\bF(N)_{G^p \times M}$ is the relative Deligne complex with respect to the fibration $\pi : G^p \times M \to G^p \times pt$.
\end{dfn}

Note that the complex of simplicial sheaves $\bF(N)_{G^{\bullet} \times M}$ can also be given by the quotient $\bF(N)_{G^{\bullet} \times M} = \F(N)_{G^{\bullet} \times M} / F^1\!\F(N)_{G^{\bullet} \times M}$. 

As is clear, if $G$ is a discrete group, then $\F(N)_{G^{\bullet} \times M} = \bF(N)_{G^{\bullet} \times M}$.

\begin{dfn}
Let $G$ be a Lie group acting on a smooth manifold $M$. We define the \textit{$G$-equivariant smooth Deligne cohomology group} of $M$ to be the hypercohomology group $H^m(G^{\bullet} \times M, \bF(N)_{G^{\bullet} \times M})$ of the complex of simplicial sheaves $\bF(N)_{G^{\bullet} \times M}$ on $G^{\bullet} \times M$.
\end{dfn}

From now on, we omit the subscripts of $\F(N)_{G^i \times M}$, $\F(N)_{G^{\bullet} \times M}$, etc. We also write $H^m(G^{\bullet} \times M, \bF(0)) = H^m(G^{\bullet} \times M, \u{\T})$.

\begin{rem}
For a finite group $G$, the hypercohomology $H^m(G^{\bullet} \times M, \F(N)) = H^m(G^{\bullet} \times M, \bF(N))$ is introduced in the work of Lupercio and Uribe \cite{L-U} as the \textit{Deligne cohomology group for the orbifold $M/G$}.
\end{rem}

By definition, the hypercohomology $H^m(G^{\bullet} \times M, \bF(N))$ is given in the following way. Let $I^{*, *, *}$ be an injective resolution of $\bF(N)$, that is, $I^{i, *, *}$ is an injective resolution of the complex of sheaves $\bF(N)$ on $G^i \times M$:
$$
\xymatrix{
\vdots &
\vdots &
&
\vdots &
\vdots \\
I^{i, 1, 0} \ar[r]^{\til{d}} \ar[u]^{\delta} &
I^{i, 1, 1} \ar[r]^{\til{d}} \ar[u]^{\delta} &
\cdots \ar[r]^-{\til{d}} &
I^{i, 1, N} \ar[r]^{\til{d}} \ar[u]^{\delta} &
I^{i, 1, N+1} \ar[r]^{\til{d}} \ar[u]^{\delta} &
\cdots \\
I^{i, 0, 0} \ar[r]^{\til{d}} \ar[u]^{\delta} &
I^{i, 0, 1} \ar[r]^{\til{d}} \ar[u]^{\delta} &
\cdots \ar[r]^-{\til{d}} &
I^{i, 0, N} \ar[r]^{\til{d}} \ar[u]^{\delta} &
I^{i, 0, N+1} \ar[r]^{\til{d}} \ar[u]^{\delta} &
\cdots \\
\u{\T} \ar[r]^{\til{d}} \ar[u] &
\u{A}^1_{rel} \ar[r]^{\til{d}} \ar[u] &
\cdots \ar[r]^-{\til{d}} &
\u{A}^N_{rel} \ar[r]^{\til{d}} \ar[u] &
0 \ar[r]^{\til{d}} \ar[u] &
\cdots.
}
$$
We define a triple complex $(K^{i, j, k}, \d, \delta, \til{d})$ by
\begin{eqnarray}
K^{i, j, k} = \Gamma(G^i \times M, I^{i, j, k}).
\label{triple_complex:resolution}
\end{eqnarray}
On $\oplus_{m = i+j+k} K^{i, j, k}$, the total coboundary operator is defined by $D = \d + (-1)^i \delta + (-1)^{i+j} \til{d}$ on the component $K^{i, j, k}$. The cohomology of this total complex is $H^m(G^{\bullet} \times M, \bF(N))$.

\subsection{The \Cech cohomology description}

It is possible to compute the equivariant smooth Deligne cohomology by means of a \Cech cohomology. For this purpose, we introduce the notion of an open cover of the simplicial manifold $G^{\bullet} \times M$.

\begin{dfn}[\cite{Bry3,Bry2}]
We define an open cover of the simplicial manifold $G^{\bullet} \times M$ by a family of open covers $\U^{\bullet} = \{ \U^{(p)} \}_{p \ge 0}$ such that:

\begin{enumerate}
\item
$\U^{(p)} = \{ U^{(p)}_{\alpha^{(p)}} \}_{\alpha^{(p)} \in \AA^{(p)} }$ is an open cover of $G^p \times M$;

\item
the index set $\AA^{(p)}$ forms a simplicial set $\AA^{\bullet} = \{ \AA^{(p)} \}_{p \ge 0}$; and

\item
we have $\d_i( U^{(p+1)}_{\alpha^{(p+1)}} ) \subset U^{(p)}_{\d_i( \alpha^{(p+1)} )}$ and $s_i( U^{(p)}_{\alpha^{(p)}} ) \subset U^{(p+1)}_{s_i( \alpha^{(p)} )}$.
\end{enumerate}
\end{dfn}

For an open cover $\U^{\bullet} = \{ \U^{(p)} \}_{p \ge 0}$ of $G^{\bullet} \times M$, we define a triple complex $(\check{K}^{i, j, k}, \d, \delta, \til{d})$ by
\begin{eqnarray}
\check{K}^{i, j, k} = 
\prod_{U^{(i)}_{\alpha^{(i)}_0}, \ldots, U^{(i)}_{\alpha^{(i)}_j}}
\Gamma(U^{(i)}_{\alpha^{(i)}_0 \ldots \alpha^{(i)}_j}, \bF(N)^{[k]}),
\label{triple_complex:Cech}
\end{eqnarray}
where we write $\bF(N)^{[k]}$ for the sheaf located at degree $k$ in the complex $\bF(N)$. The coboundary operator $\d : \check{K}^{i, j, k} \to \check{K}^{i+1, j, k}$ is given by $\d = \sum_{l=0}^{i+1}(-1)^l\d^*_l$, the $\delta : \check{K}^{i, j, k} \to \check{K}^{i, j+1, k}$ is the \Cech coboundary operator, and $\til{d} : \check{K}^{i, j, k} \to \check{K}^{i, j, k+1}$ is induced by the coboundary operator $\til{d} : \bF(N)^{[k]} \to \bF(N)^{[k+1]}$. From the triple complex, we obtain the total complex by putting $\check{C}^m(\U^{\bullet}, \bF(N)) = \oplus_{m=i+j+k}\check{K}^{i, j, k}$, where the total coboundary operator is defined by $D = \d + (-1)^i \delta + (-1)^{i+j}\til{d}$ on the component $\check{K}^{i, j, k}$. We denote by $\check{H}^m(\U^{\bullet}, \bF(N))$ the cohomology of the total complex.

As in the case of ordinary sheaf cohomology, there exists a canonical homomorphism $\check{H}^m(\U^{\bullet}, \bF(N)) \to H^m(G^{\bullet} \times M, \bF(N))$. This induces an isomorphism 
$$
\lim_{\to} \check{H}^m(\U^{\bullet}, \bF(N)) \longrightarrow 
H^m(G^{\bullet} \times M, \bF(N)),
$$
where the direct limit is taken over the ordered set of open covers of $G^{\bullet} \times M$. 

When each open cover $\U^{(p)}$ is a \textit{good cover} \cite{B-T,Bry1} of $G^p \times M$, we call $\U^{\bullet}$ a \textit{good cover} of $G^{\bullet} \times M$.

\begin{lem}
If $\U^{\bullet}$ is a good cover of $G^{\bullet} \times M$, then there exists a natural isomorphism $\check{H}^m(\U^{\bullet}, \bF(N)) \cong H^m(G^{\bullet} \times M, \bF(N))$.
\end{lem}

\begin{proof}
Because $\U^{(p)}$ is a good cover of $G^p \times M$ for each $p$, we have an isomorphism $\check{H}^q(\U^{(p)}, \bF(N)) \cong H^q(G^p \times M, \bF(N))$. Using the spectral sequence associated with the filtration $F^pK = \oplus_{i \ge p}K^{i, *, *}$ and that associated with $F^p\check{K} = \oplus_{i \ge p}\check{K}^{i, *, *}$, we obtain the isomorphism in the lemma.
\end{proof}

We give an example of an open cover \cite{Bry2}. Let $\U = \{ U_\alpha \}_{\alpha \in \AA}$ be an open cover of $M$. For $p = 0, 1, \ldots$ we put $\AA^{(p)} = \AA^{p+1}$. We define face maps $\d_i : \AA^{(p+1)} \to \AA^{(p)}$ and degeneracy maps $s_i : \AA^{(p)} \to \AA^{(p+1)}$ by
\begin{eqnarray*}
\d_i(\alpha_0, \ldots, \alpha_{p+1}) & = &
(\alpha_0, \ldots, \widehat{\alpha_i}, \ldots, \alpha_{p+1}), \\
s_i(\alpha_0, \ldots, \alpha_{p}) & = &
(\alpha_0, \ldots, \alpha_i, \alpha_i, \ldots, \alpha_p).
\end{eqnarray*}
The open cover $\U^{(p)} = \{ U^{(p)}_{\alpha^{(p)}} \}_{\alpha^{(p)} \in \AA^{(p)}}$ is inductively defined by
$$
U^{(p)}_{\alpha^{(p)}} = 
\bigcap_{i=0}^p \d_i^{-1}( U^{(p-1)}_{\d_i(\alpha^{(p)})}).
$$


\section{Properties of equivariant smooth Deligne cohomology groups}
\label{sec:properties_EDC}

We study some general properties of the equivariant smooth Deligne cohomology group $H^m(G^{\bullet} \times M, \bF(N))$ in relations with other cohomology groups. Throughout this section, $G$ denotes a Lie group acting on a smooth manifold $M$.

\subsection{Relation to smooth Deligne cohomology}

\begin{prop} \label{prop:EDC_to_DC}
There exists a spectral sequence converging to a graded quotient of $H^m(G^{\bullet} \times M, \bF(N))$ with its $E_2$-term given by
$$
E^{p, q}_2 = H^p(H^q(G^* \times M, \bF(N)), \d),
$$
the $p$th cohomology group of the complex
$$
H^q(M, \F(N)) \stackrel{\d}{\longrightarrow}
H^q(G \times M, \bF(N)) \stackrel{\d}{\longrightarrow}
H^q(G^2 \times M, \bF(N)) \stackrel{\d}{\longrightarrow} \cdots. 
$$
\end{prop}

\begin{proof}
This spectral sequence is given by the filtration $F^pK = \oplus_{i \ge p} K^{i, *, *}$ of the triple complex (\ref{triple_complex:resolution}). Then the $E_1$-terms are $E^{p, q}_1 = H^q(G^p \times M, \bF(N))$, and the differential $d_1 : E^{p, q}_1 \to E^{p+1, q}_1$ is $d_1 = \d = \sum_{i=0}^{p+1}(-1)^i\d_i^*$. Because $K^{i, j, k}$ is zero unless $i, j, k \ge 0$, the spectral sequence converges to the graded quotient of $H^{p+q}(G^{\bullet} \times M, \bF(N))$ with respect to the filtration.
\end{proof}

\begin{cor} \label{cor:G_trivial}
If $G = \{ e \}$, then $H^m(G^{\bullet} \times M, \bF(N)) \cong H^m(M, \F(N))$.
\end{cor}

\begin{proof}
The natural identification $G^p \times M = M$ implies that $E^{p, q}_1 = E^{0, q}_1$ for all $p$ and $q$. It is easy to see that $d_1 = 0$ if $p$ is even, and $d_1 = \id$ if $p$ is odd. Thus, the spectral sequence degenerates at $E_2$, and gives the result.
\end{proof}

Note that we can identify the $E_2$-term $E_2^{p, 0}$ of the spectral sequence in Proposition \ref{prop:EDC_to_DC} as the differentiable cohomology \cite{H-M} of $G$ with coefficients the $G$-module $H^0(M, \T)$:
$$
E_2^{p, 0} = H^p_{group}(G, H^0(M, \T)).
$$
A condition on $G$ allows us to express the other $E_2$-terms in a similar fashion.

\begin{lem}
If $G$ is discrete, then we can identify the $E_2$-term $E_2^{p, q}$ of the spectral sequence in Proposition \ref{prop:EDC_to_DC} as the group cohomology of degree $p$ with coefficients the smooth Deligne cohomology $H^q(M, \F(N))$ regarded as a $G$-module:
$$
E_2^{p, q} = H^p_{group}(G, H^q(M, \F(N))).
$$
In particular, we have
$$
E_2^{0, q} = H^q(M, \F(N))^G =
\{ c \in H^q(M, \F(N)) |\ g^*c = c \ \mbox{for all} \ g \in G \}.
$$
\end{lem}

\begin{proof}
Since $G$ is discrete, we have $\bF(N) = \F(N)$ and $E_1^{p, q} = H^q(G^p \times M, \F(N)) = C^p_{group}(G, H^q(M, \F(N)))$. The identification of the differential $d_1$ with the coboundary operator on the group cochains establishes the lemma.
\end{proof}

\begin{rem}
In general, we have $E_2^{0, q} \subset H^q(M, \F(N))^G$ for any Lie group $G$. However, when $G$ is not discrete, it happens that $E_2^{0, q}$ does not coincide with $H^q(M, \F(N))^G$. In fact, when $G = SU(2)$ acts on $M = SU(2)$ by the left translation, a $G$-invariant integral 3-form on $M$ is a class in $H^2(M, \F(2))^G$ which does not belong to $E_2^{0, 2}$.
\end{rem}

\subsection{Relation to equivariant cohomology}

We study a relation between the equivariant smooth Deligne cohomology and the ordinary equivariant cohomology via $H^m(G^{\bullet} \times M, \u{\T}) \cong H^m(G^{\bullet} \times M, \bF(0))$.

We denote by $EG \times_G M$ the topological space obtained by the quotient of $EG \times M$ under the diagonal action of $G$, where $EG$ is the total space of the universal $G$-bundle $EG \to BG$. As in \cite{A-B}, we define the \textit{equivariant cohomology group} $H^m_G(M, \Z)$ by $H^m_G(M, \Z) = H^m(EG \times_G M, \Z)$, where the latter is the ordinary (singular) cohomology with coefficients $\Z$.

\smallskip

Notice the following short exact sequence of simplicial sheaves on $G^{\bullet} \times M$:
$$
0 \longrightarrow
\Z \longrightarrow
\u{\R} \stackrel{\exp2\pi\i}{\longrightarrow}
\u{\T} \longrightarrow 0,
$$
where $\Z$ is the constant simplicial sheaf, and $\u{\R}$ is given by the sheaf of germs of $\R$-valued functions on each $G^p \times M$. Thus, by the associated long exact sequence, we have a homomorphism 
$$
H^m(G^{\bullet} \times M, \u{\T}) \longrightarrow 
H^{m+1}(G^{\bullet} \times M, \Z).
$$

\begin{lem}[\cite{Bry2}] 
If $G$ is compact, then there exists a natural isomorphism $H^m(G^{\bullet} \times M, \u{\T}) \cong H^{m+1}_G(M, \Z)$ for $m > 0$.
\end{lem}

\begin{proof}
We have an exact sequence
$$
H^m(G^{\bullet} \times M, \u{\R}) \to
H^m(G^{\bullet} \times M, \u{\T}) \to
H^{m+1}(G^{\bullet} \times M, \Z) \to
H^{m+1}(G^{\bullet} \times M, \u{\R}).
$$
It is proved by Dupont \cite{Du} that there exists a natural isomorphism $H^m(G^{\bullet} \times M, \Z) \cong H^m(|G^{\bullet} \times M|, \Z)$. Since $|G^{\bullet} \times M| \cong EG \times_G M$, the lemma will follow from the vanishing of $H^m(G^{\bullet} \times M, \u{\R})$ for $m > 0$. 

There is a spectral sequence converging to a graded quotient of $H^{p+q}(G^{\bullet} \times M, \u{\R})$ with its $E_1$-term given by
$$
E^{p, q}_1 = H^q(G^p \times M, \u{\R}) =
\left\{
\begin{array}{cl}
C^{\infty}(G^p \times M, \R), & \quad (q = 0), \\
0, & \quad (q > 0).
\end{array}
\right.
$$
Here we used the fact that $\u{\R}$ is soft, which is a consequence of the existence a partition of unity on $G^p \times M$. This spectral sequence degenerates at $E_2$, and we have $H^p(G^{\bullet} \times M, \u{\R}) \cong E^{p, 0}_2$.

Now, since $G$ is compact, we have an invariant measure $dg$ on $G$. We suppose that the measure is normalized. For $f \in E^{p, 0}_1$, we define $\bar{f} \in E^{p-1, 0}_1$ by
$$
\bar{f}(g_1, \ldots, g_{p-1}, x) =
\int_{g \in G}f(g, g_1, \ldots, g_{p-1}, x) dg.
$$
If $\d f = 0$, then $f = \d \bar{f}$. Hence we have $E^{p, 0}_2 = H^p(G^{\bullet} \times M, \u{\R}) = 0$ for $p > 0$.
\end{proof}

As a generalization of (\ref{exact_seq:Deligne_to_integral:relative}), we consider a short exact sequence of complexes of simplicial sheaves on $G^{\bullet} \times M$:
\begin{eqnarray}
0 \to 
\{ 0 \to \u{A}^1_{rel} \to \cdots \to \u{A}^N_{rel} \} \to
\bF(N) \to
\{ \u{\T} \to 0 \to \cdots \to 0 \} \to 0. 
\label{exact_seq:Deligne_to_integral:simplicial}
\end{eqnarray}
Composing the induced homomorphism 
$$
H^m(G^{\bullet} \times M, \bF(N)) \longrightarrow 
H^m(G^{\bullet} \times M, \u{\T})$$ 
with the homomorphism $H^m(G^{\bullet} \times M, \u{\T}) \to H^{m+1}(G^{\bullet} \times M, \Z) \cong H^{m+1}_G(M, \Z)$, we obtain a homomorphism to the equivariant cohomology group
$$
H^m(G^{\bullet} \times M, \bF(N)) \longrightarrow H^{m+1}_G(M, \Z).
$$

\begin{lem}
If $G$ is compact, then $H^m(G^{\bullet} \times M, 0 \to \u{A}^1_{rel} \to \cdots \to \u{A}^N_{rel})$ is
$$
\left\{
\begin{array}{cl}
A^1(M)^G_{cl}, & m = 1, \\
A^m(M)^G_{cl} / dA^{m-1}(M)^G, & 1 < m < N, \\
A^N(M)^G / dA^{N-1}(M)^G, & m = N, \\
0, & \mbox{otherwise},
\end{array}
\right.
$$
where $A^p(M)^G$ is the group of $G$-invariant $p$-forms on $M$, and $A^p(M)^G_{cl}$ is the group of $G$-invariant closed $p$-forms on $M$.
\end{lem}

\begin{proof}
Because the sheaf $\u{A}^p_{rel}$ is soft on each $G^i \times M$, the group of our interest is computed as the cohomology of the double complex $(L^{*, *}, \d, \til{d})$, where 
$$
L^{i, j} = 
\left\{
\begin{array}{cl}
A^j(G^i \times M)_{rel}, & \quad 1 \le j \le N, \\
0, & \quad \mbox{otherwise}.
\end{array}
\right.
$$
We compute this cohomology by the spectral sequence associated with the filtration $F^pL = \oplus_{j \ge p} L^{*, j}$. By using the invariant measure on $G$, we obtain 
$$
E^{p, q}_1 = H^q(L^{*, p}, \d) =
\left\{
\begin{array}{cl}
A^p(M)^G, & \quad q = 0, \\
0, & \quad q > 0.
\end{array}
\right.
$$
Thus, $E^{p, q}_2 = 0$ for $q > 0$, and $E^{p, 0}_2$ is the $p$th cohomology of the complex
$$
0 \longrightarrow
A^1(M)^G \stackrel{d}{\longrightarrow}
A^2(M)^G \stackrel{d}{\longrightarrow}
\cdots \stackrel{d}{\longrightarrow}
A^N(M)^G \longrightarrow 
0 \longrightarrow \cdots.
$$
The spectral sequence degenerates at $E_2$, and we proved the lemma.
\end{proof}

\begin{prop} \label{prop:EDC_to_EC}
Suppose that $G$ is compact and $N$ is a positive integer. 

(a) The group $H^N(G^{\bullet} \times M, \bF(N))$ fits into the exact sequence
$$
\xymatrix@C=10pt@R=10pt{
H^N_G(M, \Z) \ar[d] \\
A^N(M)^G / d A^{N-1}(M)^G \ar[r] &
H^N(G^{\bullet} \times M, \bF(N)) \ar[r] & 
H^{N+1}_G(M, \Z) \ar[r] & 0.
} 
$$

(b) If $m \ge N$, then $H^m(G^{\bullet} \times M, \bF(N)) \cong H^{m+1}_G(M, \Z)$.
\end{prop}

\begin{proof}
These are straight consequences of the long exact sequence associated with (\ref{exact_seq:Deligne_to_integral:simplicial}) and the lemmas above.
\end{proof}

\subsection{Relation to invariant differential forms}

We have the following short exact sequence of complexes of simplicial sheaves, which is a generalization of (\ref{exact_seq:Deligne_to_forms:relative}):
\begin{eqnarray}
\xymatrix@C=7pt{
0 \ar[r] &
\{ \u{\T} \to \u{A}^1_{rel} \to \cdots \to \u{A}^N_{rel, cl} \} \ar[r] &
\bF(N) \ar[r]^-{d} &
\{ 0 \to \cdots \to 0 \to \u{A}^{N+1}_{rel, cl} \} \ar[r] & 0.
} \label{exact_seq:Deligne_to_forms:simplicial}
\end{eqnarray}
By the relative \Poincare lemma \cite{Bry1}, we have a quasi-isomorphism
$$
\{ \pi^{-1}\u{\T} \to 0 \to \cdots \to 0 \} \to
\{ \u{\T} \to \u{A}^1_{rel} \to \cdots \to \u{A}^N_{rel,cl} \},
$$
where the simplicial sheaf $\pi^{-1}\u{\T}$ on $G^{\bullet} \times M$ is the inverse image of the simplicial sheaf $\u{\T}$ on $G^{\bullet} \times pt$ under the projection $\pi : G^{\bullet} \times M \to G^{\bullet} \times pt$.

\begin{prop} \label{prop:EDC_to_invariant_forms}
Suppose that $N$ is a positive integer.

(a) If $0 \le m < N$, then $H^m(G^{\bullet} \times M, \bF(N)) \cong H^m(G^{\bullet} \times M, \pi^{-1}\u{\T})$.

(b) The group $H^N(G^{\bullet} \times M, \bF(N))$ fits into the exact sequence
$$
\xymatrix@C=10pt@R=10pt{
0 \ar[r] &
H^N(G^{\bullet} \times M, \pi^{-1}\u{\T}) \ar[r] &
H^N(G^{\bullet} \times M, \bF(N)) \ar[r]^-d & 
A^{N+1}(M)^G_{cl} \ar[d] \\ 
& & &
H^{N+1}(G^{\bullet} \times M, \pi^{-1}\u{\T}).
} 
$$
\end{prop}

\begin{proof}
By the relative version of the \Poincare lemma, we take $(\u{A}^{*+N+1}_{rel}, d)$ as a resolution of $\u{A}^{N+1}_{rel, cl}$. By a direct computation, we obtain
$$
H^m(G^{\bullet} \times M, 0 \to \cdots \to 0 \to \u{A}^{N+1}_{rel,cl}) =
\left\{
\begin{array}{cl}
0, & (0 \le m < N), \\
A^{N+1}(M)^G_{cl}, & (m = N). \\
\end{array}
\right.
$$
Now the long exact sequence associated with (\ref{exact_seq:Deligne_to_forms:simplicial}) leads to the proposition.
\end{proof}

\subsection{A map to the equivariant de Rham cohomology}

For the smooth Deligne cohomology group $H^N(M, \F(N))$, there exists a homomorphism to the de Rham cohomology group $H^{N+1}_{DR}(M)$. This is given by taking the de Rham cohomology class of the image under the homomorphism $d : H^N(M, \F(N)) \to A^{N+1}_{cl}(N)$ induced by (\ref{exact_seq:Deligne_to_forms:manifold}). As an equivariant analogy, we here construct a homomorphism from the equivariant smooth Deligne cohomology $H^N(M, \bF(N))$ to the equivariant de Rham cohomology. 

As is known, we can formulate the ``equivariant de Rham cohomology'' by various models. The model that we employ is the following simplicial model.

\begin{dfn}
We define the \textit{equivariant de Rham cohomology group} by the total cohomology of the double complex  $(A^j(G^i \times M), \d, d)$.
\end{dfn}

\begin{rem}
By the extended de Rham theorem proved by Bott, Shulman and Stasheff \cite{B-S-S} (see also \cite{Du}), the equivariant de Rham cohomology introduced above is isomorphic to the equivariant cohomology with coefficients $\R$:
$$
H^m(A^*(G^* \times M), \d, d) \cong H^m_G(M, \R).
$$
\end{rem}

The method of constructing a map to the equivariant de Rham cohomology group is to extend the image of the homomorphism $d : H^N(G^{\bullet} \times M, \bF(N)) \to A^{N+1}(M)^G_{cl}$ induced by (\ref{exact_seq:Deligne_to_forms:simplicial}).  As the part of the ``equivariant extension'', we use the image of the homomorphism 
$$
\beta : 
H^N(G^{\bullet} \times M, \bF(N)) \longrightarrow
H^{N+1}(G^{\bullet} \times M, F^1\!\F(N))
$$ 
induced by the short exact sequence of complexes of simplicial sheaves
\begin{eqnarray}
0 \longrightarrow 
F^1\!\F(N) \longrightarrow 
\F(N) \longrightarrow 
\bF(N) \longrightarrow 0.
\label{exact_seq:1st_filtration}
\end{eqnarray}
Since $G^p \times M$ is assumed to admit a partition of unity for each $p$, it is easy to see that $H^m(G^{\bullet} \times M, F^1\!\F(N))$ is the $m$th cohomology of the double complex $(L^{i, j}, \d, \til{d})$ given by
\begin{eqnarray}
L^{i, j} = 
\left\{
\begin{array}{cl}
F^1\!A^j(G^i \times M), & \quad (1 \le j \le N), \\
0, & \quad \mbox{otherwise}.
\end{array}
\right. \label{double_complex:F1FN}
\end{eqnarray}

\begin{prop} \label{prop:EDC_to_EDRC}
There exists a homomorphism
$$
\Upsilon : H^N(G^{\bullet} \times M, \bF(N)) \to
H^{N+1}(A^*(G^* \times M), \d , d).
$$
\end{prop}

\begin{proof}
Consider a homomorphism of short exact sequences of complexes of simplicial sheaves on $G^{\bullet} \times M$ given in Figure \ref{fig:hom_short_exact_seq}.
\begin{figure}[htbp]
$$
\xymatrix{
0 \ar[d] & 
0 \ar[d] \\
F^1\!\F(N) \ar[r]^-d \ar[d] &
\{ 0 \to \cdots \to 0 \to F^1\!\u{A}^{N+1} \} \ar[d] \\
\F(N) \ar[r]^-d \ar[d] &
\{ 0 \to \cdots \to 0 \to \u{A}^{N+1} \} \ar[d] \\
\bF(N) \ar[r]^-d \ar[d] &
\{ 0 \to \cdots \to 0 \to \u{A}^{N+1}_{rel} \} \ar[d] \\
0 &
0.
}
$$
\caption{The homomorphism of short exact sequences}
\label{fig:hom_short_exact_seq}
\end{figure}
The homomorphism yields the following commutative diagram:
\begin{eqnarray}
\xymatrix{
H^N(G^{\bullet} \times M, \bF(N)) \ar[r]^d \ar[d]^{\beta} &
A^{N+1}(M)^G \ar[d]^{\d} \\
H^{N+1}(G^{\bullet} \times M, F^1\!\F(N)) \ar[r]^-d &
H^{N+1}(G^{\bullet} \times M, 0 \to \cdots \to 0 \to \u{A}^{N+1}_{rel}).
}
\end{eqnarray}
The cohomology $H^{N+1}(G^{\bullet} \times M, 0 \to \cdots \to 0 \to \u{A}^{N+1}_{rel})$ is 
$$
Ker\{ \d : F^1\!A^{N+1}(G \times M) \to F^1\!A^{N+1}(G^2 \times M) \}.
$$
Now, for a class $c \in H^N(G^{\bullet} \times M, \bF(N))$, we denote by $F^{(0)}$ the $G$-invariant $(N+1)$-form on $M$ obtained by applying $d : H^N(G^{\bullet} \times M, \bF(N)) \to A^{N+1}(M)^G$. Proposition \ref{prop:EDC_to_invariant_forms} implies that $d F^{(0)} = 0$. Because $H^{N+1}(G^{\bullet} \times M, F^1\!\F(N))$ is the total cohomology of (\ref{double_complex:F1FN}), we represent $\beta(c) \in H^{N+1}(G^{\bullet} \times M, F^1\!\F(N))$ by
$$
(F^{(1)}, \ldots, F^{(N)}) \in
F^1\!A^{N}(G \times M) \oplus \cdots \oplus F^1\!A^{1}(G^N \times M)
$$
such that $\d F^{(i-1)} + (-1)^i dF^{(i)} = 0$ for $i = 2, \ldots, N$ and $\d F^{(N)} = 0$. By the commutative diagram, we have $\d F^{(0)} = d F^{(1)}$. Hence we obtain an $(N+1)$-cocycle of the double complex $(A^*(G^* \times M), \d , d)$:
$$
(F^{(0)}, F^{(1)}, \ldots, F^{(N)}, 0) \in
A^{N+1}(M) \oplus A^N(G \times M) \oplus \cdots \oplus A^0(G^{N+1} \times M).
$$
It is straightforward to check that the equivariant de Rham cohomology class represented by the $(N+1)$-cocycle above is independent of the choice of the representative of $\beta(c)$. Thus, we obtain the homomorphism $\Upsilon$ by putting $\Upsilon(c) = [(F^{(0)}, F^{(1)}, \ldots, F^{(N)}, 0)]$.
\end{proof}


\section{Equivariant circle bundles and gerbes}
\label{sec:equiv_T_bundle_and_gerbe}

In this section, we classify equivariant principal $\T$-bundles with connection  and equivariant gerbes with connection, by using equivariant smooth Deligne cohomology groups. The essential idea here owes also to Brylinski's paper \cite{Bry2}.

We denote by $G$ a Lie group acting on a smooth manifold $M$.

\subsection{Equivariant circle bundle with connection}

As usual, a $G$-equivariant principal $\T$-bundle over $M$ is defined to be a principal $\T$-bundle $P \to M$ together with a lift of the $G$-action on $M$ to that on $P$ by bundle isomorphisms. Following \cite{Bry2}, we first reformulate this definition by using the simplicial manifold $G^{\bullet} \times M$.

For a principal $\T$-bundle $P \to M$, we define a principal $\T$-bundle $\d P \to G \times M$ by $\d P = \d_0^* P \otimes \d_1^* P^{\otimes -1}$. Similarly, for a $\T$-bundle $Q \to G \times M$, we define a $\T$-bundle $\d Q \to G^2 \times M$ by $\d Q = \d_0^*Q \otimes \d_1^*Q^{\otimes -1} \otimes \d_2^*Q$. By the relation (\ref{rel1}), the $\T$-bundle $\d \d P \to G^2 \times M$ is canonically isomorphic to the trivial $\T$-bundle $(G^2 \times M) \times \T \to G^2 \times M$.
\begin{lem} \label{lem:lift_of_action}
For a principal $\T$-bundle $P$ over $M$, the following notions are equivalent:
\begin{enumerate}
\item
a lift of the $G$-action on $M$ to that on $P$ by bundle isomorphisms;

\item 
a section $\sigma \in \Gamma(G \times M, \d P)$ such that $\d \sigma = 1$ on $G^2 \times M$, where we put $\d \sigma = \d_0^*\sigma \otimes \d_1^*\sigma^{\otimes -1} \otimes \d_2^*\sigma$.
\end{enumerate}
\end{lem}

\begin{proof}
Let $\pi : P \to M$ be a principal $\T$-bundle with a lift of the $G$-action on $M$ to that on $P$ by bundle isomorphisms. By definition, we have
\begin{eqnarray*}
\d_0^* P & = & \{ (g, x, p) \in G \times M \times P |\  x = \pi(p) \}, \\
\d_1^* P & = & \{ (g, x, p) \in G \times M \times P |\  gx = \pi(p) \}.
\end{eqnarray*}
Using the lift of the $G$-action, we define a bundle isomorphism $\varrho : \d_0^* P \to \d_1^* P$ by $\varrho(g, x, p) = (g, x, g p)$. The inverse $\varrho^{-1}$ gives rise to the section $\sigma$. Since $(g_1 g_2)p = g_1 (g_2 p)$ holds for $g_1, g_2 \in G$ and $p \in P$, we have $\d_2^*\varrho \circ \d_0^*\varrho = \d_1^*\varrho$, so that $\d \sigma = 1$. By performing this construction conversely, we can obtain a lift of the $G$-action from a section $\sigma : G \times M \to \d P$ such that $\d \sigma = 1$.
\end{proof}

\begin{thm}[Brylinski \cite{Bry2}] \label{thm:classify_equiv_bundle}
The isomorphism classes of $G$-equivariant $\T$-bundles over $M$ are classified by $H^1(G^{\bullet} \times M, \u{\T})$.
\end{thm}

The proof of this theorem can be seen in that of Theorem \ref{thm:classify_equiv_bundle_connection}.

\medskip

Next, we reformulate the notion of an invariant connection on an equivariant principal $\T$-bundle. We denote the Lie algebra of $G$ by $\g$, and the natural contraction by $\langle \ | \ \rangle : \g \otimes \mathrm{Hom}(\g, \i\R) \to \i\R$.

\begin{dfn} \label{dfn:moment_map}
Let $P \to M$ be a $G$-equivariant $\T$-bundle, and $\theta$ a connection on $P$ (which is not necessarily $G$-invariant). We define a function $\mu : M \to \mathrm{Hom}(\g, \i\R)$ by $\langle X | \mu(x) \rangle = \theta(p; X^*)$, where $p \in P_x$ is a point on the fiber of $x$, and $X^* \in T_pP$ is the tangent vector generated by the infinitesimal action of $X \in \g$ on $P$.
\end{dfn}

Since the left action of $G$ on $P$ commutes with the right action of $\T$ on $P$, the map $\mu$ is well-defined. We call the above map $\mu$ the \textit{moment} \cite{B-V}.

For a connection $\theta$ on a principal $\T$-bundle $P$ on $M$, we denote by $\d \theta$ the induced connection $\d_0^* \theta \otimes \d_1^* \theta^{\otimes -1}$ on $\d P = \d_0^*P \otimes \d_1^*P^{\otimes -1}$.

\begin{lem}
Let $P \to M$ be a $G$-equivariant $\T$-bundle, $\sigma : G \times M \to \d P$ the induced section, and $\theta$ a connection on $P$ which is not necessarily $G$-invariant. For a tangent vector $gX \oplus V \in T_gG \oplus T_xM$, the value of the 1-form $\sigma^*(\d \theta)$ is
\begin{eqnarray}
\left( \sigma^*(\d \theta) \right) 
((g, x); gX \oplus V) = 
(\theta - g^*\theta)(x; V) - \langle X | \mu(gx) \rangle.
\label{formula:expression_connection}
\end{eqnarray}
\end{lem}

\begin{proof}
By the help of the bundle map $\varrho : \d_0^*P \to \d_1^*P$ used in the proof of Lemma \ref{lem:lift_of_action}, we evaluate the value of the 1-form at the tangent vector to give
\begin{eqnarray*}
&   &
\left( \sigma^*(\d\theta) \right) 
((g, x); gX \oplus V) \\
&   & \quad \quad =
(\d_0^*\theta - \varrho^*\d_1^*\theta) 
((g, x, p); gX \oplus V \oplus \til{V}) \\
&   & \quad \quad =
\d_0^*\theta((g, x, p); gX \oplus V \oplus \til{V}) - 
\d_1^*\theta((g, x, gp); gX \oplus V \oplus (gX^* + g\til{V})) \\
&   & \quad \quad =
(\theta - g^*\theta)(p; \til{V}) - (g^*\theta)(p; X^*),
\end{eqnarray*}
where $\til{V} \in T_pP$ is a lift of the tangent vector $V \in T_xM$. Rewriting the last expression, we obtain the result.
\end{proof}

We denote by $[\omega]_{rel} \in A^q(G^p \times M)_{rel} = A^q(G^p \times M) / F^1\!A^q(G^p \times M)$ the relative differential form represented by $\omega \in A^q(G^p \times M)$.

\begin{lem} \label{lem:invariant_connection}
A connection $\theta$ on a $G$-equivariant principal $\T$-bundle $P \to M$ is $G$-invariant if and only if $[\sigma^*(\d \theta)]_{rel} = 0$ in $A^1(G \times M)_{rel}$.
\end{lem}

\begin{proof}
When $\theta$ is a $G$-invariant connection, the expression (\ref{formula:expression_connection}) implies that the 1-form $\sigma^*(\d\theta)$ belongs to $F^1\!A^1(G \times M)$. Thus, as an element in $A^1(G \times M)_{rel}$, the 1-form is equal to zero. The converse is apparent.
\end{proof}

\begin{thm} \label{thm:classify_equiv_bundle_connection}
The isomorphism classes of $G$-equivariant $\T$-bundles over $M$ with $G$-invariant connection are classified by $H^1(G^{\bullet} \times M, \bF(1))$.
\end{thm}

\begin{proof}
First, for an equivariant $\T$-bundle with connection $(P, \theta)$ over $M$, we define a cohomology class in $H^1(G^{\bullet} \times M, \bF(1))$. Let $\U^{\bullet} = \{ \U^{(p)} \}$ be a sufficiently fine open cover of $G^{\bullet} \times M$ so that we can take local sections $s_{\alpha} : U^{(0)}_{\alpha} \to P|_{U^{(0)}_{\alpha}}$. We define a cochain
\begin{eqnarray}
\left(
\begin{array}{c}
f_{\alpha_0 \alpha_1} \\
\theta^1_{\alpha}, \quad \ g_{\alpha}
\end{array}
\right)
\in \check{C}^1(\U^{\bullet}, \bF(1)) = 
\begin{array}{c}
\check{K}^{0, 1, 0} \oplus \\
\check{K}^{0, 0, 1} \oplus 
\check{K}^{1, 0, 0}
\end{array}
\label{cochain:degree_one}
\end{eqnarray}
by setting
\begin{eqnarray}
s_{\alpha_1} & = & 
s_{\alpha_0} \ f_{\alpha_0 \alpha_1}, \\
\theta^1_{\alpha} & = &
\frac{1}{2\pi\i}s_{\alpha}^* \theta, \\
(\d s)_\alpha & = & \sigma|_{ U^{(1)}_{\alpha} } \ g_{\alpha}.
\label{cocycle:100}
\end{eqnarray}
In the last line, we put $(\d s)_\alpha = \d_0^* s_{\d_0(\alpha)} \otimes \d_1^* s_{\d_1(\alpha)}^{\otimes -1}$. The cochain is closed in $\check{C}^1(\U^{\bullet}, \bF(1))$, and defines a class in $H^1(G^{\bullet} \times M, \bF(1))$. It is straight to verify that the class is independent of the choice of the local sections and of the open cover. We can also verify that the cohomology class is identical for the other equivariant $\T$-bundle with connection isomorphic to $(P, \theta)$.

Secondly, we show the injectivity of this assignment of a cohomology class. If the cocycle is a coboundary, then we obtain a global $G$-invariant section $s$ of $P$ such that $s^*\theta = 0$. Trivializing $P$ by this section, we see that the lift of the $G$-action on $M$ to that on $P \simeq M \times \T$ is trivial. Hence, $(P, \theta)$ is isomorphic to the trivial equivariant $\T$-bundle with the trivial connection. 

Lastly, we show the surjectivity. As a part of a cocycle in $\check{C}^1(\U^{\bullet}, \bF(1))$, we have a cocycle $(f_{\alpha_0 \alpha_1}, \theta^1_\alpha) \in \check{C}^1(\U^{(0)}, \F(1))$. By the standard method, we can construct a $\T$-bundle $P$ with a connection $\theta$ from the cocycle on $M$. By using $g_\alpha$ and (\ref{cocycle:100}), we can construct a global section $\sigma : G \times M \to P$ which makes $(P, \theta)$ into $G$-equivariant.
\end{proof}

We can directly generalize the proof above to obtain the classification of equivariant principal $\T$-bundles with \textit{flat} connection.

\begin{cor}
The isomorphism classes of $G$-equivariant $\T$-bundles over $M$ with $G$-invariant flat connection are classified by $H^1(G^{\bullet} \times M, \bF(N))$, where $N$ is an integer such that $N > 1$.
\end{cor}

Theorem \ref{thm:classify_equiv_bundle_connection} allows us to have a generalization of some results of Hattori and Yoshida \cite{H-Y}. 

\begin{cor} \label{cor:obstruction_moduli_bundle_connection}
Let $G$ be a Lie group acting on a smooth manifold $M$, and $P$ a principal $\T$-bundle over $M$ equipped with a connection $\theta$.

(a) There exists two obstruction classes for $(P, \theta)$ to being $G$-equivariant. The first obstruction class belongs to $H^1(G \times M, \bF(1))$, and the second obstruction class to $H^2_{group}(G, H^0(M, \T))$.

(b) Suppose that $(P, \theta)$ admits a lift of the $G$-action on $M$ by bundle isomorphisms preserving the connection. Such liftings are, up to automorphisms of $(P, \theta)$, in one to one correspondence with $H^1_{group}(G, H^0(M, \T))$.
\end{cor}

\begin{proof}
By the spectral sequence in Proposition \ref{prop:EDC_to_DC}, we obtain an exact sequence
$$
0 \longrightarrow
E^{1, 0}_2 \longrightarrow
H^1(G^{\bullet} \times M, \bF(1)) \longrightarrow
E^{0, 1}_2 \longrightarrow
E^{2, 0}_2.
$$
Recall that $(P, \theta)$ is classified by $E^{0, 1}_1 = H^1(M, \F(1))$. Because $E^{0,1}_2 = Ker \{ d_1 : E^{0, 1}_1 \to E^{1, 1}_1 \}$, the first obstruction belongs to $E^{1, 1}_1 = H^1(G \times M, \bF(1))$. Clearly, the second obstruction belongs to $E^{2, 0}_2 = H^2_{group}(G, H^0(M, \T))$. Thus (a) is proved. The (b) is also clear by $E^{1, 0}_2 = H^1_{group}(G, H^0(M, \T))$.
\end{proof}

Note that when $G$ is a finite group, Corollary \ref{cor:obstruction_moduli_bundle_connection} (b) was shown by Sharpe as the classification of orbifold group actions on a $U(1)$ gauge field \cite{Sha}.

\medskip

It would be worth while deriving the moment $\mu : M \to \mathrm{Hom}(\g, \i\R)$ in the context of the equivariant smooth Deligne cohomology group. Let $\g$ be the Lie algebra of $G$, $\g^*$ the dual space of $\g$, and $\langle \ | \ \rangle : \g \otimes \g^* \to \R$ the natural contraction. By the (co)adjoint, the Lie group $G$ acts on an element $f \in \g^*$ by $\langle X | \Ad_{g} f \rangle = \langle \Ad_{g^{-1}} X| f \rangle$.

\begin{lem} 
There exists an isomorphism
\begin{eqnarray}
H^2(G^{\bullet} \times M, F^1\!\F(1)) \cong
\left\{ f : M \to \g^* |\ 
f(gx) = \Ad_g f (x) \ \mbox{for all} \ g \in G \right\}.
\label{iso_H2F1F1}
\end{eqnarray}
\end{lem}

\begin{proof}
Since $F^1\!\u{A}^1$ is soft on $G^p \times M$ for each $p$, we have
$$
H^2(G^{\bullet} \times M, F^1\!\F(1)) =
Ker\{ \d : F^1\!A^1(G \times M) \to F^1\!A^1(G^2 \times M)\}.
$$
Let $\alpha$ be an element in $F^1\!A^1(G \times M)$. Note that for any tangent vector $V \in T_xM$ we have $\alpha((g, x); V) = 0$. By a computation, we can see that the cocycle condition $\d \alpha = 0$ is equivalent to the following conditions:
\begin{eqnarray*}
\alpha((g_2, x); g_2X) & = & \alpha((g_1g_2, x); g_1g_2X), \\
\alpha((g_1, g_2x); g_1X) & = & \alpha((g_1, g_2x); g_1Xg_2),
\end{eqnarray*}
where a tangent vector at $g \in G$ is expressed as $gX \in T_gG$ by an element $X \in T_eG = \g$. Thus, the isomorphism (\ref{iso_H2F1F1}) is induced by the assignment to $\alpha$ of the map $f : M \to \g^*$ defined by $\langle X | f(x) \rangle = \alpha((e, x); X)$. Clearly, the inverse homomorphism is given by the assignment to $f : M \to \g^*$ of the 1-form $\alpha$ defined by $\alpha((g, x); gX \oplus V) = \langle X | f(x) \rangle$.
\end{proof}

Recall that the exact sequence of complexes of simplicial sheaves
$$
0 \longrightarrow 
F^1\!\F(1) \longrightarrow
\F(1) \longrightarrow
\bF(1) \longrightarrow
0
$$
induces a homomorphism
$$
\beta : H^1(G^{\bullet} \times M, \bF(1)) \longrightarrow
H^2(G^{\bullet} \times M, F^1\!\F(1)).
$$

\begin{lem} \label{lem:beta_c:degree1}
Let $(P, \theta)$ be a $G$-equivariant $\T$-bundle with connection over $M$, and $c \in H^1(G^{\bullet} \times M, \bF(1))$ the cohomology class that classifies $(P, \theta)$. The image $\beta(c)$ is identified with the map $\frac{-1}{2\pi\i} \mu : M \to \g^*$ under (\ref{iso_H2F1F1}).
\end{lem}

\begin{proof}
We use the same notations as in the proof of Theorem \ref{thm:classify_equiv_bundle_connection}. We represent the class $c \in H^1(G^{\bullet} \times M, \bF(1))$ by the \Cech cocycle (\ref{cochain:degree_one}) in $\check{C}^1(\U^{\bullet}, \bF(1))$. Regarding it as a cochain in $\check{C}^1(\U^{\bullet}, \F(1))$, we compute the coboundary. Then the non-trivial components are given by
$$
\d \theta^1_\alpha - \frac{1}{2\pi\i}d\log g_\alpha 
\in A^1(U^{(1)}_\alpha),
$$
where $\d \theta^1_\alpha = \d_0^* \theta^1_{\d_0(\alpha)} - \d_1^* \theta^1_{\d_1(\alpha)}$. By the definition of $\theta^1_\alpha$ and $g_\alpha$, we obtain
$$
\d \theta^1_\alpha - \frac{1}{2\pi\i}d\log g_\alpha =
\frac{1}{2\pi\i} 
\sigma^*(\d\theta)|_{U^{(1)}_\alpha} 
\in F^1\!A^1(U^{(1)}_\alpha).
$$
Now (\ref{formula:expression_connection}) and (\ref{iso_H2F1F1}) complete the proof.
\end{proof}

We remark that the computations in the proof above is essentially the same as that performed in \cite{Bry2} to obtain the equivariant Chern class of $(P, \theta)$.

\subsection{Equivariant gerbe with connection}

The notion of equivariant gerbes (with connection) \cite{Bry2} is defined in a fashion similar to the simplicial formulation of equivariant principal $\T$-bundles (with connection). To save pages, we follow terminologies in \cite{Bry1,Bry-M}, and drop including the definition of gerbe itself.

As in the case of principal $\T$-bundles, for a gerbe $\sC$ over $M$ we define a gerbe $\d \sC$ over $G \times M$ by setting $\d \sC = \d_0^*\sC \otimes \d_1^*\sC^{\otimes -1}$. Similar notations will be used in the sequel.

\begin{dfn}[\cite{Bry2}] \label{dfn:equiv_gerbe}
A \textit{$G$-equivariant gerbe} is defined to be a gerbe $\sC$ over $M$ together with the following data:
\begin{enumerate}
\item[(1)] 
a global object $R \in \Gamma(G \times M, \d\sC)$;

\item[(2)] 
an isomorphism $\psi : \d R \to \mathbf{1}$ of global objects in the trivial gerbe over $G^2 \times M$ which satisfies $\d\psi = 1$ on $G^3 \times M$.
\end{enumerate}
\end{dfn}

In the definition above, we regard the trivial gerbe over $G^2 \times M$ as the sheaf of categories of principal $\T$-bundles. So the global object $\mathbf{1}$ means the trivial $\T$-bundle $(G^2 \times M) \times \T \to G^2 \times M$.

\begin{dfn}[\cite{Bry2}] \label{dfn:equiv_trivial_gerbe}
A $G$-equivariant gerbe $(\sC, R, \psi)$ is said to be \textit{trivial} if we have a global object $P \in \Gamma(M, \sC)$ and an isomorphism $\eta : R \to \d P$ on $G \times M$ such that $\d\eta = \psi$ on $G^2 \times M$. The data $(P, \eta)$ is called an \textit{equivariant object} of $(\sC, R, \psi)$.
\end{dfn}

\begin{thm}[Brylinski \cite{Bry2}]  \label{thm:classify_equiv_gerbe}
The isomorphism classes of $G$-equivariant \\ gerbes over $M$ are classified by $H^2(G^{\bullet} \times M, \u{\T})$.
\end{thm}

The proof of this theorem can be seen in that of Theorem \ref{thm:classify_equiv_gerbe_connection}.

\begin{dfn}[\cite{Bry2}] \label{dfn:equiv_gerbe_connection}
Let $(\sC, R, \psi)$ be a $G$-equivariant gerbe over $M$. 

(a) A \textit{$G$-invariant connective structure} on $(\sC, R, \psi)$ is defined to be a connective structure $\Co$ on $\sC$ with an element $D_{rel} \in \Gamma(G \times M, (\d\Co)_{rel}(R))$ such that the isomorphism $\psi$ carries $\d D_{rel}$ to the relatively trivial connection on $\mathbf{1}$.

(b) A \textit{$G$-invariant curving} for the $G$-invariant connective structure $(\Co, D_{rel})$ is defined to be a curving $K$ for $\Co$ such that $(\d K)_{rel}(D_{rel}) = 0$.
\end{dfn}

Here we add a few explanations to the definition above. A connective structure $\Co$ on $\sC$ induces a connective structure $\d \Co = \d_0^*\Co \otimes \d_1^*\Co^{\otimes -1}$ on $\d \sC$, and a curving $K$ for $\Co$ induces a curving $\d K = \d_0^* K \otimes \d_1^*K^{\otimes -1}$ for $\d \Co$. From the connective structure, we obtain the \textit{relative connective structure} by setting 
$$
(\d \Co)_{rel}(P) =
(\d \Co)(P) / F^1\!\u{A}^1
$$
for a local object $P \in \Gamma(U, \d \sC)$ given on an open set $U \subset G \times M$. We denote by $[\nabla]_{rel} \in \Gamma(U, (\d \Co)_{rel}(P)) = \Gamma(U, (\d \Co)(P)) / F^1\!A^1(U)$ the element represented by $\nabla \in \Gamma(U, \d \Co(P))$. For $(\d \Co)_{rel}$ the \textit{relative curving}
$$
(\d K)_{rel} : 
\Gamma(U, (\d \Co)_{rel}(P)) \to A^2(U)_{rel}
$$
is defined by $(\d K)_{rel}([\nabla]_{rel}) = [(\d K)(\nabla)]_{rel}$.

When a $G$-equivariant gerbe is equipped with a $G$-invariant connective structure and a $G$-invariant curving, we simply call it a $G$-equivariant gerbe with connection.

\begin{dfn}
A $G$-equivariant gerbe with connection $(\sC, \Co, K, R, D_{rel}, \psi)$ is said to be trivial if we have a global object $P \in \Gamma(M, \sC)$, an element $\nabla \in \Gamma(M, \Co(P))$ and an isomorphism $\eta : (R, D_{rel}) \to (\d P, [\d \nabla]_{rel})$ such that $K(\nabla) = 0$ and $\d \eta = \psi$.
\end{dfn}

\begin{thm} \label{thm:classify_equiv_gerbe_connection}
The isomorphism classes of $G$-equivariant gerbes with connection over $M$ are classified by $H^2(G^{\bullet} \times M, \bF(2))$.
\end{thm}

\begin{proof}
First, we assign a \Cech cocycle to an equivariant gerbe with connection $(\sC, \Co, K, R, D_{rel}, \psi)$. Let $\U^{\bullet} = \{ \U^{(p)}\}$ be a sufficiently fine open cover of $G^{\bullet} \times M$. For each open set $U^{(0)}_{\alpha}$, we take an object $P_{\alpha} \in \Gamma(U^{(0)}_{\alpha}, \sC)$ and an element $\nabla_{\alpha} \in \Gamma(U^{(0)}_{\alpha}, \Co(P_{\alpha}))$. We also take isomorphisms $u_{\alpha_0 \alpha_1} : P_{\alpha_1} \to P_{\alpha_0}$ on $U^{(0)}_{\alpha_0 \alpha_1}$ and $v_{\alpha} : R \to (\d P)_\alpha = \d_0^*P_{\d_0(\alpha)} \otimes \d_1^*P^{\otimes -1}_{\d_0(\alpha)}$ on $U^{(1)}_{\alpha}$. We use below the following notations:
\begin{eqnarray*}
(\d u)_{\alpha_0 \alpha_1} & = & 
\d_0^* u_{\d_0(\alpha_0) \d_0(\alpha_1)} \otimes 
\d_1^* u_{\d_1(\alpha_0) \d_1(\alpha_1)}^{\otimes -1}, \\
(\d \nabla)_{\alpha} & = &
\d_0^* \nabla_{\d_0(\alpha)} \otimes 
\d_1^* \nabla_{\d_1(\alpha)}^{\otimes -1}, \\
(\d v)_\alpha & = &
\d_0^* v_{\d_0(\alpha)} \otimes
\d_1^* v_{\d_1(\alpha)}^{\otimes -1} \otimes
\d_2^* v_{\d_2(\alpha)}.
\end{eqnarray*}
Then we define a \Cech cochain
\begin{eqnarray}
\left(
\begin{array}{c}
f_{\alpha_0 \alpha_1 \alpha_2} \\
\theta^1_{\alpha_0 \alpha_1},
g_{\alpha_0 \alpha_1} \\
\theta^2_{\alpha}, \quad \
\omega^1_{\alpha}, \quad \
h_{\alpha} 
\end{array}
\right)
\in \check{C}^2(\U^{\bullet}, \bF(2)) = 
\begin{array}{c}
\check{K}^{0, 2, 0} \oplus \\
\check{K}^{0, 1, 1} \oplus 
\check{K}^{1, 1, 0} \oplus \\
\check{K}^{0, 0, 2} \oplus 
\check{K}^{1, 0, 1} \oplus
\check{K}^{2, 0, 0}
\end{array}
\label{cochain:degree2}
\end{eqnarray}
by setting
\begin{eqnarray*}
f_{\alpha_0 \alpha_1 \alpha_2}
& = &
u_{\alpha_0 \alpha_1} \circ 
u_{\alpha_1 \alpha_2} \circ 
u_{\alpha_2 \alpha_0}, \\
\theta^1_{\alpha_0 \alpha_1}
& = &
(u_{\alpha_0 \alpha_1})_* \nabla_{\alpha_1} - \nabla_{\alpha_0}, \\
\theta^2_{\alpha} & = & K(\nabla_{\alpha}), \\
g_{\alpha_0 \alpha_1}
& = &
v_{\alpha_0}^{-1} \circ (\d u)_{\alpha_0 \alpha_1} \circ v_{\alpha_1}, \\
\omega^1_{\alpha}
& = &
[(\d \nabla)_\alpha]_{rel} - (v_\alpha)_* D_{rel}, \\
h_{\alpha} 
& = &
\psi \circ (\d v)_\alpha^{-1}.
\end{eqnarray*}
By computations, we can verify that the cochain is indeed a cocycle. We can also verify that the cohomology class represented by the cocycle is independent of all the choices made, and is identical for the other equivariant gerbe with connection isomorphic to $(\sC, \Co, K, R, D_{rel}, \psi)$.

Secondly, we show that the assignment of a cohomology class is injective. We suppose that the cocycle assigned to $(\sC, \Co, K, R, D_{rel}, \psi)$ is a coboundary. Without loss of generality, we can assume that the cocycle is zero. Then we have $(f_{\alpha_0 \alpha_1 \alpha_2}, \theta^1_{\alpha_0 \alpha_1}, \theta^2_\alpha) = 0 \in \check{C}^2(\U^{(0)}, \F(2))$. This implies that $(\sC, \Co, K)$ is trivial: we have a global object $P' \in \Gamma(M, \sC)$ and an element $\nabla' \in \Gamma(M, \Co(P'))$ such that $K(\nabla') = 0$. Similarly, because $g_{\alpha_0 \alpha_1} = 1$, there is an isomorphism $\eta : R \to \d P'$. It is direct to see that $\omega^1_\alpha = 0$ and $h_\alpha = 1$ lead to the rest of the conditions for the equivariant gerbe with connection to be trivial.

Finally, we show the surjectivity. Suppose that a class in $H^2(G^{\bullet} \times M, \bF(2))$ is given. We represent it by a \Cech cocycle as in (\ref{cochain:degree2}). Then we can construct a gerbe with connection $(\sC, \Co, K)$ from the cocycle condition $D^{(0)} (f, \theta^1, \theta^2) = 0$ in $\check{C}^3(\U^{(0)}, \F(2))$, where $D^{(0)}$ is the coboundary operator on $\check{C}^2(\U^{(0)}, \F(2))$. (See \cite{Bry1} for the detailed construction.) The condition $\d (f, \theta^1, \theta^2) - D^{(1)} (g, \omega^1) = 0$ in $\check{C}^2(\U^{(1)}, \bF(2))$ implies that we have a global object $R \in \Gamma(G \times M, \d \sC)$ and an element $D_{rel} \in \Gamma(G \times M, (\d \Co)_{rel}(R))$ such that $(\d K)_{rel}(D_{rel}) = 0$. By using the condition $\d(g, \omega^1) + D^{(2)}(h) = 0$ in $\check{C}^1(\U^{(2)}, \bF(2))$, we obtain a global isomorphism $\psi : \d R \to \mathbf{1}$ which carries $\d D_{rel}$ to the relatively trivial connection on $\mathbf{1}$. This isomorphism satisfies $\d \psi = 1$ by the rest of the condition $\d (h) = 1$ in $\check{C}^0(\U^{(3)}, \bF(2))$.
\end{proof}

A similar proof establishes the following corollary.

\begin{cor}
The isomorphism classes of $G$-equivariant gerbes with flat connection over $M$ are classified by $H^2(G^{\bullet} \times M, \bF(N))$, where $N$ is an integer such that $N > 2$.
\end{cor}

By using a result in Section \ref{sec:properties_EDC}, we obtain the gerbe version of Corollary \ref{cor:obstruction_moduli_bundle_connection}.

\begin{cor} \label{cor:obstruction_moduli_gerbe_connection}
Let $G$ be a Lie group acting on a smooth manifold $M$, and $(\sC, \Co, K)$ a gerbe with connection over $M$. 

(a) There are three obstruction classes for $(\sC, \Co, K)$ to being $G$-equivariant. The first obstruction class belongs to $H^2(G \times M, \bF(2))$. The second and third obstruction classes are represented by cohomology classes in $H^1(G^2 \times M, \bF(2))$ and $H^3_{group}(G, H^0(M, \T))$ respectively.

(b) Suppose that we can make $(\sC, \Co, K)$ into a $G$-equivariant gerbe with connection over $M$. The isomorphism classes of such $G$-equivariant gerbes with connection are in one to one correspondence with a subgroup $F^1\!H^2$ contained in $H^2(G^{\bullet} \times M, \bF(2))$ which fits into the following exact sequence
\begin{eqnarray*}
\xymatrix@C=15pt@R=3pt{
H^0(H^1(G^* \times M, \bF(2)), \d) \ar[r] &
H^2_{group}(G, H^0(M, \T)) \ar[r] &
F^1\!H^2 \ar[r] & \ar@{} \\
H^1(H^1(G^* \times M, \bF(2)), \d) \ar[r] &
H^3_{group}(G, H^0(M, \T)).
}
\end{eqnarray*}
\end{cor}

\begin{proof}
By the spectral sequence in Proposition \ref{prop:EDC_to_DC}, we have two exact sequences
\begin{eqnarray*}
\xymatrix@C=15pt@R=3pt{
0 \ar[r] & 
F^1\!H^2(G^{\bullet} \times M, \bF(2)) \ar[r] &
H^2(G^{\bullet} \times M, \bF(2)) \ar[r] &
E^{0, 2}_{\infty} \ar[r] &
0, \\
0 \ar[r] & 
E^{2, 0}_{\infty} \ar[r] &
F^1\!H^2(G^{\bullet} \times M, \bF(2)) \ar[r] &
E^{1, 1}_{\infty} \ar[r] &
0. 
}
\end{eqnarray*}
Recall that $(\sC, \Co, K)$ is classified by $E^{0, 2}_1 = H^2(M, \F(2))$. Thus, by means of the first exact sequence above, it suffices to understand when an element in $E^{0, 2}_1$ survives into $E^{0, 2}_{\infty}$. Note that $E^{0, 2}_{\infty} = E^{0, 2}_4 = Ker \{ d_3 : E^{0, 2}_3 \to E^{3, 0}_3 \}$ and $E^{0, 2}_3 = Ker \{ d_2 : E^{0, 2}_2 \to E^{2, 1}_2 \}$. Hence we have three obstruction classes belonging to the following groups
\begin{eqnarray*}
E^{1, 2}_1 & = & H^2(G \times M, \bF(2)), \\
E^{2, 1}_2 & = &
Ker \{ d_1 : E^{2, 1}_1 \to E^{3, 1}_1 \}/
Im \{ d_1 : E^{1, 1}_1 \to E^{2, 1}_1 \}, \\
E^{3, 0}_3 & = & E^{3, 0}_2 / Im \{ d_2 : E^{1, 1}_2 \to E^{3, 0}_2 \}.
\end{eqnarray*}
Since $E^{2, 1}_1 = H^1(G^2 \times M, \bF(2))$ and $E^{3, 0}_2 = H^3_{group}(G, H^0(M, \T))$, we have (a). The first exact sequence above also implies that there is a one to one correspondence between $F^1\!H^2(G^{\bullet} \times M, \bF(2))$ and the set of isomorphism classes of $G$-equivariant gerbes with connection whose underlying gerbes with connection are $(\sC, \Co, K)$. By expressing the second short exact sequence in terms of $E_2$-terms, we obtain (b).
\end{proof}

\begin{rem}
Let $G = SU(2)^\delta$ be the group $SU(2)$ with the discrete topology. We consider the action of $G$ on $M = SU(2) = S^3$ by the left translation. By computations, we have $E^{1, 1}_1 = E^{2, 1}_1 = 0$ for the spectral sequence in Proposition \ref{prop:EDC_to_DC}. Thus, the obstruction class for a gerbe with connection over $S^3$ which is classified by an $SU(2)^\delta$-invariant class in $H^2(S^3, \F(2)) \cong A^3(S^3)_0$ to being $SU(2)^\delta$-equivariant belongs to $E^{3, 0}_2 = H^3_{group}(SU(2)^{\delta}, \T)$. A detailed description of the obstruction class is given by Brylinski \cite{Bry1}.
\end{rem}

When $G$ is a finite group, the choice of the ways of making a gerbe with connection over $M$ into $G$-equivariant was intensively studied by Sharpe  \cite{Sha} as the classification of orbifold group actions on B-fields. In \cite{L-U}, Lupercio and Uribe performed the classification by using the Deligne cohomology group for the orbifold $M/G$, which coincides with the equivariant smooth Deligne cohomology group in the case of $G$ finite.

As Corollary \ref{cor:obstruction_moduli_gerbe_connection} indicates, we can ``twist'' a $G$-equivariant gerbe with connection by a group 2-cocycle with coefficients in the $G$-module $H^0(M, \T)$: let $(\sC, \Co, K, R, D_{rel}, \psi)$ be a $G$-equivariant gerbe with connection over $M$, and $\gamma : G \times G \times M \to \T$ a smooth function which defines a group 2-cocycle $\gamma \in Z^2_{group}(G, H^0(M, \T))$. If we replace the isomorphism $\psi$ of principal $\T$-bundles by $\psi \gamma$, then $(\sC, \Co, K, R, D_{rel}, \psi \gamma)$ is also an equivariant gerbe with connection.

\medskip

By using the other result in Section \ref{sec:properties_EDC}, we have the next corollary.

\begin{cor}
Let $G$ be a Lie group acting on a smooth manifold $M$. If $G$ is compact, then any $G$-equivariant gerbe over $M$ admits a $G$-invariant connective structure and a $G$-invariant curving.
\end{cor}

\begin{proof}
Recall the homomorphism $H^2(G^{\bullet} \times M, \bF(2)) \to H^2(G^{\bullet} \times M, \u{\T})$. Since $G$ is assumed to be compact, the homomorphism is surjective by Proposition \ref{prop:EDC_to_EC} (a). Now the corollary follows from Theorem \ref{thm:classify_equiv_gerbe} and Theorem \ref{thm:classify_equiv_gerbe_connection}.
\end{proof}

For a gerbe with connection $(\sC, \Co, K)$ over $M$, we have a closed 3-form $\Omega \in A^3(M)_{cl}$ called the \textit{3-curvature}. It is proved by Brylinski \cite{Bry1} that a closed 3-form $\Omega \in A^3(M)_{cl}$ is the 3-curvature of a gerbe with connection if and only if it is integral, namely $\Omega$ belongs to the kernel of the map $A^3(M)_{cl} \to H^3(M, \T)$ induced from (\ref{exact_seq:Deligne_to_forms:manifold}). 

\begin{cor}
Let $G$ be a Lie group acting on a smooth manifold $M$.

(a) For a $G$-equivariant gerbe with connection $(\sC, \Co, K, R, D_{rel}, \psi)$ over $M$, the 3-curvature $\Omega \in A^3(M)_{cl}$ of $(\sC, \Co, K)$ is $G$-invariant.

(b) A $G$-invariant closed 3-form $\Omega \in A^3(M)_{cl}^G$ is the 3-curvature of a $G$-equivariant gerbe with connection if and only if $\Omega$ belongs to the kernel of the homomorphism $A^3(M)^G_{cl} \to H^3(G^{\bullet} \times M, \pi^{-1}\u{\T})$ induced from (\ref{exact_seq:Deligne_to_forms:simplicial}).
\end{cor}

\begin{proof}
By Proposition \ref{prop:EDC_to_invariant_forms} (b), we have an exact sequence

$$
H^2(G^{\bullet} \times M, \bF(2)) \stackrel{d}{\longrightarrow}
A^3(M)^G_{cl} \longrightarrow
H^3(G^{\bullet} \times M, \pi^{-1}\u{\T}).
$$
Since the 3-curvature of $(\sC, \Co, K)$ is obtained as the image under $d$ of the class in $H^2(G^{\bullet} \times M, \bF(2))$ that classifies $(\sC, \Co, K, R, D_{rel}, \psi)$, the exact sequence above directly leads to the corollary.
\end{proof}

As is shown in Proposition \ref{prop:EDC_to_EDRC}, the image of the homomorphism $\beta : H^2(G^{\bullet} \times M, \bF(2) \to H^3(G^{\bullet} \times M, F^1\!\F(2))$ provides the ``equivariant extension part'' for the 3-curvature of an equivariant gerbe with connection. By Lemma \ref{lem:beta_c:degree1}, we can think of it as a ``moment'' associated with an equivariant gerbe with connection. The image of $\beta$ is used in \cite{Go} to study a relationship between $G$-equivariant gerbes with connection over $M$ and gerbes with connection over the quotient space $M/G$.


\bigskip

\begin{acknowledgments}
The essential ideas in Brylinski's paper \cite{Bry2} influenced my work a great deal, and I would like to express my gratitude to him. I am indebted to my supervisors, Professors T. Kohno and M. Furuta for helpful discussions and good advice. I am grateful to Y. Hashimoto for reading earlier drafts and for valuable discussions. I would also like to thank I. Satake for useful suggestions and discussions. 
This research is supported by Research Fellowship of the 
Japan Society for the Promotion of Science for Young Scientists.

\end{acknowledgments}



\begin{flushleft}
Graduate school of Mathematical Sciences, University of Tokyo, \\
Komaba 3-8-1, Meguro-Ku, Tokyo, 153-8914 Japan. \\
e-mail: kgomi@ms.u-tokyo.ac.jp
\end{flushleft}

\end{document}